\documentclass[a4paper]{amsart}
\usepackage{delimset, amssymb, enumitem, caption, mathtools}
\usepackage[margin=30mm]{geometry}
\usepackage{microtype}
\setlist{
itemsep=4pt, 
topsep=2pt, 
leftmargin=17pt, 
listparindent=11pt}
\usepackage[pagebackref, ocgcolorlinks]{hyperref}
\usepackage{xcolor}
\definecolor{BIT}{cmyk}{1, 0, 1, 0}
\AtBeginDocument{\hypersetup{citecolor=BIT}}

\usepackage{aliascnt}
\newtheorem{theorem}{Theorem}[section]
\newaliascnt{lemma}{theorem}

\aliascntresetthe{lemma}
\newaliascnt{corollary}{theorem}
\newtheorem{corollary}[corollary]{Corollary}
\aliascntresetthe{corollary}
\newaliascnt{proposition}{theorem}
\newtheorem{proposition}[proposition]{Proposition}
\aliascntresetthe{proposition}
\numberwithin{equation}{section}
\theoremstyle{remark}

\usepackage[capitalize]{cleveref}
\crefname{ineq}{Ineq.}{Ineqs.}
\Crefname{ineq}{Inequality}{Inequalities}
\creflabelformat{ineq}{~\upshape(#2#1#3)}

\makeatletter
\newcommand\creflabel[2][\@currentcounter]{%
 \crefalias{\@currentcounter}{#1}\label{#2}}
\makeatother

\allowdisplaybreaks

\usepackage[numbers, sort&compress, nonamebreak, merge, elide, longnamesfirst]{natbib}

\usepackage{url}

\usepackage{tikz}
\usetikzlibrary{decorations.pathreplacing, calc, positioning}
\tikzset{
edge/.style={semithick},
ball/.style={shape=circle, minimum size=1mm, ball color=black, inner sep=0.5},
ellipsis/.style={shape=circle, fill, inner sep=.5}}
\def\r{1}
\def\eps{.2}

\author[D.Q.B.~Tang]{Davion Q.B. Tang}
\address[Davion Q.B. Tang]{School of Mathematics and Statistics, Beijing Institute of Technology, Beijing 102400, P. R. China}
\email{davion@bit.edu.cn}

\author[D.G.L.~Wang]{David G.L. Wang}
\address[David G.L. Wang]{School of Mathematics and Statistics \& MIIT Key Laboratory of Mathematical Theory and Computation in Information Security, Beijing Institute of Technology, Beijing 102400, P. R. China}
\email{glw@bit.edu.cn}

\thanks{D.Q.B. Tang: School of Mathematics and Statistics, Beijing Institute of Technology, Beijing 102400, P. R. China; davion@bit.edu.cn. D.G.L. Wang: School of Mathematics and Statistics \& MIIT Key Laboratory of Mathematical Theory and Computation in Information Security, Beijing Institute of Technology, Beijing 102400, P. R. China; glw@bit.edu.cn (corresponding author). D.G.L. Wang is supported by the General Program of the National Natural Science Foundation of China (Grant No.~12171034).}


\title[Positive $e_I$-expansions of graphs]{Positive $e_I$-expansions for KPKP graphs, twinned lollipops, and kayak paddles}

\begin{document}

\bibliographystyle{abbrvnat}

\begin{abstract}
We derive explicit positive $e_I$-expansions for the chromatic symmetric functions of several graph families. First, we obtain a uniform formula for KPKP graphs that specializes to known formulas for lollipops, KPK graphs, KKP graphs, and PKP graphs. We then give positive $e_I$-expansions for twinned paths and cycles and apply the KPKP formula to twinned lollipops. Finally, we establish a positive $e_I$-expansion for kayak paddles and specialize it to infinity graphs. These formulas refine ordinary $e$-positivity by retaining composition-indexed coefficients and provide independent proofs complementary to existing recurrence, unit-interval, and noncommutative methods.
\end{abstract}
\maketitle

\section{Introduction}
This paper continues the study of \citet{WZ25} on the $e$-positivity of graphs. In his seminal paper, \citet{Sta95} introduced the chromatic symmetric function $X_G$ of a graph $G$ and posed the problem of characterizing $e$-positive graphs. Recently, \citet{Hik24X} proved that every unit interval graph is $e$-positive, thereby confirming the Stanley--Stembridge conjecture \citep{SS93}; see also \citet{Hik25}.

Although \citet{Hik24X} provided a probabilistic formula for the $e$-coefficients of the chromatic quasisymmetric function of a unit interval graph, explicit positive expansions for individual graph families remain useful. Such expansions provide effective coefficient formulas, independent combinatorial proofs, and tools that apply beyond the unit-interval setting. Recent work on support restrictions and selected $e$-coefficients illustrates this continuing direction; see \citet{ST26}.

Several methods have been developed for proving $e$-positivity of a graph $G$. A direct method is to display the $e$-expansion of the chromatic symmetric function $X_G$, as has been done for paths and cycles; see \citet{Sta95,Wol98,SW16,Ell17}. \citet{GS01} introduced appendable $(e)$-positivity in the algebra $\mathrm{NCSym}$ of symmetric functions in noncommuting variables. Since $(e)$-positivity implies ordinary $e$-positivity, this method yields, among other results, the $e$-positivity of $K$-chains. Generating functions and recurrences provide another effective approach; for example, \citet{Dv18} used them to re-establish the $e$-positivity of lollipops.

A further approach uses the algebra $\mathrm{NSym}$ of noncommutative symmetric functions and elementary expansions indexed by compositions. Following the notation of \citet{GKLLRT95}, we use the symbol $I$ to denote a composition. Inspired by the noncommutative treatment of Schur positivity in \citet{TW23X}, \citet{WZ25} developed the \emph{composition method} for producing explicit $e_I$-expansions of chromatic symmetric functions. The method has since been applied to conjoined graphs, cycle-chords, clocks, spiders, and trinacria graphs; see \citet{QTW26,Wang25,CHW26,TWW24X,GWZ25X}. Its main advantage for the present paper is that it records useful information at the level of compositions before passing to the ordinary elementary basis.

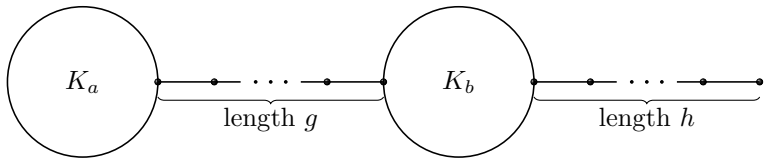
\begin{figure}[htbp]
\begin{tikzpicture}[decoration=brace]
\node (2) at (2,0) [ball]{};
\node (6) at (2.75,0) [ball]{};
\node (7) at (4.25,0) [ball]{};
\node (a) at (3.5,0) [ellipsis]{};
\node (al) at (3.5-\eps,0) [ellipsis]{};
\node (ar) at (3.5+\eps,0) [ellipsis]{};
\node (3) at (5,0) [ball]{};
\node (4) at (7,0) [ball]{};
\node (5) at (10,0) [ball]{};
\node (6) at (7.75,0) [ball]{};
\node (7) at (9.25,0) [ball]{};
\node (b) at (8.5,0) [ellipsis]{};
\node (bl) at (8.5-\eps,0) [ellipsis]{};
\node (br) at (8.5+\eps,0) [ellipsis]{};
\draw[edge] (0, 0) arc [start angle=180, end angle=540, x radius=\r, y radius=\r];
\draw[edge] (3) arc [start angle=180, end angle=540, x radius=\r, y radius=\r];
\draw[edge] (2, 0) -- ($(al) - (\eps, 0)$);
\draw[edge] ($(ar) + (\eps, 0)$) -- (5, 0);
\draw[edge] (7, 0) -- ($(bl) - (\eps, 0)$);
\draw[edge] ($(br) + (\eps, 0)$) -- (10, 0);
\node at (1,0) {$K_a$};
\node at (6,0) {$K_b$};
\draw [decorate] (5,-.2) -- (2,-.2);
\draw [decorate] (10,-.2) -- (7,-.2);
\node at (3.5,-.5) {length $g$};
\node at (8.5,-.5) {length $h$};
\end{tikzpicture}
\caption{The KPKP graph $P^g(K_a,\,K_b^h)$.}\label{fig:KPKP}
\end{figure}

This paper has three parts. First, we derive a positive $e_I$-expansion for the chromatic symmetric function of KPKP graphs, in which the letters K and P indicate complete graphs and paths, respectively; see \cref{fig:KPKP,sec:KPKP}.
The formula specializes to the known formulas for lollipops obtained by \citet{Tom25}, for KPK graphs obtained by \citet{WZ25}, and for KKP and PKP graphs obtained by \citet{QTW26}; see \cref{thm:lollipop.melting,thm:KPK,thm:KPK:PKP}. Every KPKP graph is a $K$-chain and was therefore already known to be $e$-positive. Our result gives a different, family-specific positive $e_I$-expansion, obtained through algebraic and combinatorial arguments for functions on compositions; see \cref{thm:KPKP}.

\begin{figure}[htbp]
\begin{tikzpicture}[decoration=brace]
\node (2) at (2, 0) [ball]{};
\node (8) at (2.75, 0) [ball]{};
\node (9) at (4.25, 0) [ball]{};
\node (a) at (3.5, 0) [ellipsis]{};
\node (al) at (3.5-\eps, 0) [ellipsis]{};
\node (ar) at (3.5+\eps, 0) [ellipsis]{};
\node (3) at (5, 0) [ball]{};
\node (4) at (6, 0) [ball]{};
\node (5) at (7, 0) [ball]{};
\node (11) at (8, 0) [ball]{};
\node (b) at (9, 0) [ellipsis]{};
\node (bl) at (9-\eps, 0) [ellipsis]{};
\node (br) at (9+\eps, 0) [ellipsis]{};
\node (6) at (10, 0) [ball]{};
\node (10) at (11, 0) [ball]{};
\node (7) at (6, .6) [ball]{};
\draw[edge] (0,0) arc [start angle=180, end angle=540, x radius=\r, y radius=\r];
\draw[edge] (2, 0) -- ($(al) - (\eps, 0)$);
\draw[edge] ($(ar) + (\eps, 0)$) -- (5, 0);
\draw[edge] (6, .6) -- (5, 0);
\draw[edge] (6, .6) -- (6, 0);
\draw[edge] (6, .6) -- (7, 0);
\draw[edge] (5, 0) -- (7, 0);
\draw[edge] (7, 0) -- ($(bl) - (\eps, 0)$);
\draw[edge] ($(br) + (\eps, 0)$) -- (11, 0);
\node at (1,0) {$K_a$};
\draw [decorate] (11,-.2) -- (6,-.2);
\node at (8.5,0)[below=7pt] {length $h$};
\draw [decorate] (5,-.2) -- (2,-.2);
\node at (3.5,0)[below=7pt] {length $l-h-1$};
\end{tikzpicture}
\caption{The twinned lollipop $\mathrm{tw}_h(K_a^l)$.}
\label{fig:lollipop.tw}
\end{figure}
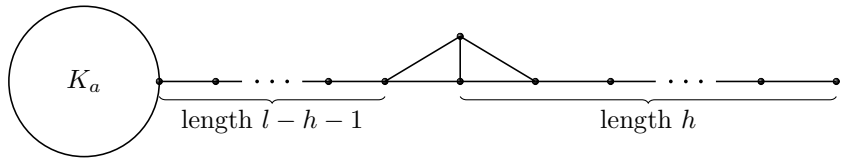

Second, we give positive $e_I$-expansions for the chromatic symmetric functions of twinned paths and twinned cycles; see \cref{thm:path.tw,thm:cycle.tw}. The $e$-positivity of these graphs was established by \citet{BCCCGKKLLS25} using generating functions and recurrences. The graph operation of twinning at a vertex was introduced by \citet{FHM19}, who conjectured that twinning preserves $e$-positivity. \citet{LLWY21} produced counterexamples to this conjecture for general graphs, whereas \citet{BCCCGKKLLS25} proved preservation for paths and cycles. We continue this study by deriving an explicit positive $e_I$-expansion for every twinned lollipop; see \cref{fig:lollipop.tw}.
Twinned lollipops form a family of unit interval graphs, so their $e$-positivity also follows from \citeauthor{Hik24X}'s theorem. The contribution here is the more explicit composition-indexed formula and an independent proof by the composition method.

\begin{figure}[htbp]
\begin{tikzpicture}[decoration=brace]
\draw (0,0) node {$C_a$};
\draw (6,0) node {$C_b$};
\node[ball] (c) at (0: \r) {};
\node (p1) at ($ (c) + (\r, 0) $) [ball] {};
\node[ellipsis] at ($ (p1) + (\r-\eps, 0) $) {};
\node[ellipsis] (e2) at ($ (p1) + (\r, 0) $) {};
\node[ellipsis] at ($ (p1) + (\r+\eps, 0) $) {};
\node (p2) at ($ (p1) + (2*\r, 0) $) [ball] {};
\node (p3) at ($ (p2) + (\r, 0) $) [ball] {};
\node[below=5pt] at (e2) {length $l$};
\draw[edge] (-75: \r) arc (-75: 75: \r);
\draw[edge] (5, 0) arc [start angle=180, end angle=105, radius=\r];
\draw[edge] (5, 0) arc [start angle=180, end angle=255, radius=\r];
\draw[edge] (c) -- ($(p1) + (\r-2*\eps, 0)$);
\draw[edge] (p3) -- (p2) -- ($(p2) - (\r-2*\eps, 0)$);
\draw [decorate] (5,-.2) -- (1,-.2);
\end{tikzpicture}
\caption{The kayak paddle $P^l(C_a,C_b)$.}
\label{fig:kayak}
\end{figure}
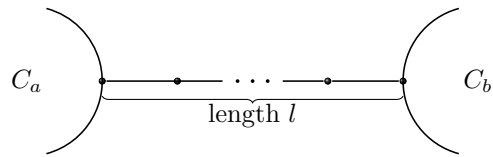

Third, we obtain an explicit positive $e_I$-expansion for the chromatic symmetric function of kayak paddles; see \cref{fig:kayak}. The name \emph{kayak paddle} was used by \citet{AWv24}, who proved the $e$-positivity of this family by an extension of \citeauthor{GS01}'s appendable $(e)$-positivity method. Our proof uses a cyclic-averaging argument and combinatorial injections.
The same family was studied by \citet{MS25} under the name \emph{dumbbell graphs}, which also appears in earlier spectral graph literature; see \citet{WHBL09}. Since the first version of this paper, \citet{TV26-HatChain,TV26-vglue} have proved the $e$-positivity of adjacent cycle--clique chains, and have established a broader theorem for graphs obtained by gluing sequences of unit interval graphs and cycles at single vertices. The latter result includes kayak paddles at the level of $e$-positivity. Our formula is complementary because it gives a concrete positive $e_I$-expansion. 

The principal contributions of this paper may be summarized as follows. First, we obtain a uniform positive $e_I$-expansion for KPKP graphs that simultaneously specializes to the known composition-indexed formulas for lollipops, KPK graphs, KKP graphs, and PKP graphs. Second, we give explicit positive $e_I$-expansions for twinned paths and twinned cycles and, by combining the KPKP formula with the twinning construction, for every twinned lollipop. Third, we derive a positive $e_I$-expansion for kayak paddles, together with an explicit specialization to infinity graphs, by using cyclic averaging and combinatorial injections. Thus, the results do more than certify ordinary $e$-positivity: they retain the finer composition-indexed coefficient data, unify several previously separate formulas, and provide independent proofs that complement the existing appendable-positivity, recurrence, and unit-interval-graph approaches.

This paper is organized as follows. \Cref{sec:pre} fixes the notation and recalls the results used later. In \cref{sec:KPKP}, we derive the positive $e_I$-expansion for KPKP graphs. In \cref{sec:twin}, we treat twinned paths and twinned cycles and then apply the KPKP formula to twinned lollipops. Finally, \cref{sec:kayak} is devoted to kayak paddles and infinity graphs.

\section{Preliminaries}\label{sec:pre}

This section fixes the notation and recalls the results on chromatic symmetric functions that will be used later. We follow the terminology of \citet{GKLLRT95} and \citet{Sta11B}. Let $n$ be a positive integer. A \emph{composition} of $n$ is a sequence of positive integers with sum~$n$, commonly denoted $I=i_1 \dotsm i_l\vDash n$, with \emph{size} $\abs{I}=n$, \emph{length} $\ell(I)=l$, and \emph{parts} $i_1,\dots,i_l$. When all parts $i_k$ have the same value $i$, we write $I=i^l$. The \emph{reversal} composition $i_l\dotsm i_1$ is written as $\overline{I}$. For convenience, we denote the $k$th last part $i_{l+1-k}$ by $i_{-k}$, and denote the composition obtained by removing the $k$th part by~$I\backslash i_k$, i.e.,
\[
I\backslash i_k
=i_1\dotsm i_{k-1}i_{k+1}\dotsm i_{-1}.
\]
Whenever an expression such as $i_{-r}$ occurs in a summation, that summation is understood to range only over compositions $I$ with $\ell(I)\ge r$. When a capital letter such as $I$ or $J$ denotes a composition, the corresponding lowercase letter with integer subscripts denotes its parts. A \emph{partition} of $n$ is a weakly decreasing sequence of positive integers with sum~$n$, denoted $\lambda=\lambda_1\lambda_2\dotsm\vdash n$.

A \emph{symmetric function} of homogeneous degree $n$ over the field $\mathbb{Q}$ of rational numbers is a bounded-degree formal power series 
\[
f(x_1, x_2, \dots)
=\sum_{
\lambda=\lambda_1 \lambda_2 \dotsm \vdash n}
c_\lambda \cdotp
x_1^{\lambda_1} 
x_2^{\lambda_2} 
\dotsm, 
\]
where $c_\lambda \in \mathbb Q$, 
that is invariant under every finite permutation of the variables $x_i$. Equivalently, it has a unique expansion
\[
f=\sum_{\lambda\vdash n}c_\lambda m_\lambda,
\]
where $m_\lambda$ is the monomial symmetric function. Let $\operatorname{Sym}^0=\mathbb{Q}$, and let $\operatorname{Sym}^n$ be the vector space of homogeneous symmetric functions of degree $n$ over~$\mathbb{Q}$. One basis of $\operatorname{Sym}^n$ consists of elementary symmetric functions $e_\lambda$ for all partitions $\lambda\vdash n$, where
\[
e_\lambda
=
e_{\lambda_1}
e_{\lambda_2}\dotsm
\quad\text{and}\quad
e_k
=
\sum_{1\le i_1<\dots<i_k} 
x_{i_1} \dotsm x_{i_k}.
\]
A symmetric function $f\in\mathrm{Sym}$ is said to be \emph{$e$-positive} if every $e_\lambda$-coefficient of $f$ is nonnegative.

For any composition $I$, there is a unique partition $\rho(I)$ which consists of the parts of $I$. This allows us to define $e_I=e_{\rho(I)}$. An \emph{$e_I$-expansion} of a symmetric function $f\in\mathrm{Sym}^n$ is an expression
\[
f=\sum_{I\vDash n} c_I e_I.
\]
The systematic use of composition-indexed $e_I$-expansions for chromatic symmetric functions was introduced by \citet{WZ25}. We call it a \emph{positive $e_I$-expansion} if $c_I\ge 0$ for all $I$. The $e$-positivity of a symmetric function can then be shown by presenting a positive $e_I$-expansion. Since $e_I$ depends only on the multiset of parts of $I$, such expansions are generally nonunique; this flexibility is central to the composition method.

\citet{Sta95} introduced the chromatic symmetric function for a simple graph $G$ as
\[
X_G
=\sum_\kappa \prod_{v \in V(G)} x_{\kappa(v)},
\]
where $\kappa$ runs over proper colorings of~$G$. For instance, the chromatic symmetric function of the complete graph~$K_n$ is $X_{K_n}=n!e_n$. It refines Birkhoff's chromatic polynomial $\chi_G(q)$ through the specialization $\chi_G(q)=X_G(1^q)$. A standard tool in the study of chromatic symmetric functions is the triple-deletion property established by \citet[Theorem 3.1, Corollaries 3.2 and 3.3]{OS14}.

\begin{proposition}[\citeauthor{OS14}]\label{prop:3del}
Let $G$ be a graph with a stable set $T$ of order $3$. Denote by $e_1$, $e_2$ and $e_3$ the edges linking the vertices in $T$. For any set $S\subseteq \{1,2,3\}$, denote by $G_S$ the graph with vertex set~$V(G)$ and edge set $E(G)\cup\{e_j\colon j\in S\}$. Then
\[
X_{G_{12}}
=X_{G_1}+X_{G_{23}}-X_{G_3}
\quad\text{and}\quad
X_{G_{123}}
=X_{G_{13}}+X_{G_{23}}-X_{G_3}.
\]
\end{proposition}

\citet[Table~1]{SW16} obtained the following formula using Stanley's generating function for Smirnov words, see also \citet[Theorem 7.2]{SW10}.
\begin{proposition}[\citeauthor{SW16}]\label{prop:path}
We have $X_{P_n}=\sum_{I\vDash n}w_Ie_I$ for any $n\ge 1$, where
\begin{equation}\label{def:w}
w_I
=i_1(i_2-1)(i_3-1)\dotsm(i_l-1)
\quad\text{if $I=i_1i_2\dotsm i_l$}.
\end{equation}
\end{proposition}
Analogously, \citet[Corollary 6.2]{Ell17X} gave a formula for the chromatic quasisymmetric function of cycles, whose $t=1$ specialization is the following. We use the convention $C_2=K_2$ whenever a formula involves a cycle of order~$2$.

\begin{proposition}[\citeauthor{Ell17X}]\label{prop:cycle}
We have $X_{C_n}=\sum_{I\vDash n}(i_1-1)w_Ie_I$ for $n\ge 2$.
\end{proposition}
Only compositions whose parts are all at least $2$ contribute to the sum for $X_{C_n}$. This set will be used repeatedly; we denote it by
\begin{equation}\label{def:Cn}
\mathcal W_n
=\{I\vDash n\colon i_1,i_2,\dots\ge 2\}.
\end{equation}

A \emph{double-rooted graph} is a triple $(G,u,v)$ that consists of a graph $G$ and two distinct vertices $u$ and $v$ of $G$, in which $u$ and $v$ are called the \emph{roots}. For any double-rooted graphs $(G_1,u_1,v_1),\dots,(G_l,u_l,v_l)$, denote by $G_1+\dots+G_l$ the graph obtained by identifying $v_i$ and $u_{i+1}$ for all $1\le i\le l-1$. It has order
\[
\abs{G_1+\dots+G_l}
=\abs{G_1}+\dots+\abs{G_l}-(l-1).
\]
In this case, each of the graphs $G_1$ and $G_l$ needs only one root. We call a graph with one root a \emph{rooted graph}. A graph $G$ is \emph{vertex-transitive} if for any two vertices $u$ and $v$ of $G$, there is an automorphism $f$ such that $f(u)=v$. We omit the root whenever the graph is vertex-transitive or the intended root is clear from context. A \emph{$K$-chain} is a graph of the form $K_{i_1}+\dots+K_{i_l}$, denoted $K_I$, where $I=i_1\dotsm i_l\in\mathcal W_n$. We use the same notation for the evident degenerate cases involving $K_1$, although those graphs are not $K$-chains under this definition.

\begin{figure}[htbp]
\begin{tikzpicture}[decoration=brace]
\node (c) at (0,0) [ellipsis]{};
\node (a1) at (-\r, 0) [ball]{};
\node (a2) at (-\r*2, 0) [ball]{};
\node (b1) at (\r, 0) [ball]{};
\node (b2) at (\r*2, 0) [ball]{};
\node[above right] at (a2) {$u$};
\node[above left] at (b2) {$v$};
\node (cr) at (0+\eps, 0) [ellipsis]{};
\node (cl) at (0-\eps, 0) [ellipsis]{};
\draw[edge] (a2) arc [start angle=0, end angle=360,
x radius=\r, y radius=\r*.75];
\draw[edge] (b2) arc [start angle=180, end angle=-180,
x radius=\r*.75, y radius=\r*.5];
\draw[edge] (a2) -- (a1) -- ($(cl) - (\eps, 0)$);
\draw[edge] (b2) -- (b1) -- ($(cr) + (\eps, 0)$);
\node at (-\r*3, 0) {$G$};
\node at (\r*2.75, 0) {$H$};
\node at (0, -.5) {length $l$};
\draw [decorate] (2, -.2) -- (-2, -.2);
\end{tikzpicture}
\caption{The $l$-conjoined graph $P^l(G,H)$.}
\label{fig:GH}
\end{figure}
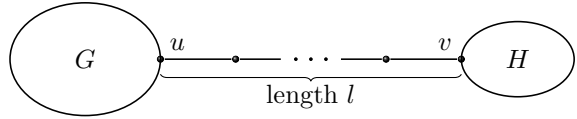

For any rooted graphs $(G,u)$ and $(H,v)$, the \emph{$l$-conjoined graph} $P^l(G,H)$ is the graph obtained by adding a path of length $l$ that links $u$ and $v$, see \cref{fig:GH}.
It has order $\abs{G}+\abs{H}+l-1$. We call $G$ and~$H$ the \emph{node graphs}. When $G=K_1$ and $H$ is vertex-transitive, we introduce the more compact notation
\[
H^l=P^l(H,K_1)
=H+P^{l+1}.
\]
In terms of the chromatic symmetric functions of graphs of the form $H^l$, \citet[Propositions~3.1 and~3.2]{QTW26} expressed the chromatic symmetric function of conjoined graphs in which one of the node graphs is a clique or a cycle.

\begin{proposition}[\citeauthor{QTW26}]\label{prop:KPG+CPG}
Let $l\ge 0$ and $a\ge 2$. For any rooted graph $H$,
\begin{align}
\label{fml:KPG}
X_{P^l(K_a,\,H)}
&=(a-1)!
\sum_{i=0}^{a-1}
(1-i)e_i
X_{H^{a+l-i-1}}
\quad\text{and}
\\
\label{fml:CPG}
X_{P^l(C_a,\,H)}
&=(a-1)X_{H^{a+l-1}}
-\sum_{i=1}^{a-2}
X_{C_{a-i}}
X_{H^{i+l-1}}.
\end{align}
\end{proposition}
The first identity also holds for $a=1$, since $P^l(K_1,H)=H^l$ and its right side reduces to~$X_{H^l}$.

\begin{figure}[htbp]
\begin{tikzpicture}[decoration=brace]
\draw (0,0) node {$K_a$};
\node[ball] (c) at (0: \r) {};
\node[ball] at (30: \r) {};
\node[ball] at (-30: \r) {};
\node (p1) at ($ (c) + (\r, 0) $) [ball] {};
\node[ellipsis] at ($ (p1) + (\r-\eps, 0) $) {};
\node[ellipsis] (e2) at ($ (p1) + (\r, 0) $) {};
\node[ellipsis] at ($ (p1) + (\r+\eps, 0) $) {};
\node (p2) at ($ (p1) + (2*\r, 0) $) [ball] {};
\draw [decorate] (4, -.2) -- (1, -.2);
\node at (2.5, -.5) {length $l$};
\draw[edge] (-75: \r) arc (-75: 75: \r);
\draw[edge] (c) -- ($(p1) + (\r-2*\eps, 0)$);
\draw[edge] (p2) -- ($(p2) - (\r-2*\eps, 0)$);

\begin{scope}[xshift=7cm]
\draw (0,0) node {$C_a$};
\node[ball] (c) at (0: \r) {};
\node[ball] at (30: \r) {};
\node[ball] at (-30: \r) {};
\node (p1) at ($ (c) + (\r, 0) $) [ball] {};
\node[ellipsis] at ($ (p1) + (\r-\eps, 0) $) {};
\node[ellipsis] (e2) at ($ (p1) + (\r, 0) $) {};
\node[ellipsis] at ($ (p1) + (\r+\eps, 0) $) {};
\node (p2) at ($ (p1) + (2*\r, 0) $) [ball] {};
\draw [decorate] (4, -.2) -- (1, -.2);
\node at (2.5, -.5) {length $l$};
\draw[edge] (-75: \r) arc (-75: 75: \r);
\draw[edge] (c) -- ($(p1) + (\r-2*\eps, 0)$);
\draw[edge] (p2) -- ($(p2) - (\r-2*\eps, 0)$);
\end{scope}
\end{tikzpicture}
\caption{The lollipop $K_a^l$ and the tadpole $C_a^l$.}
\label{fig:lollipop-tadpole}
\end{figure}
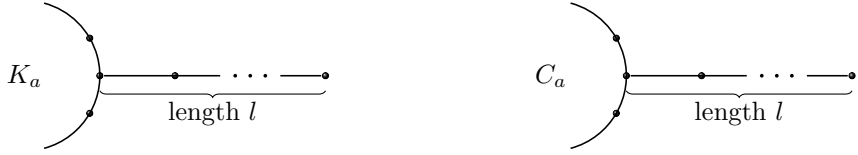

A \emph{lollipop} is a $K$-chain of the form $K_a^l$, see the left illustration in \cref{fig:lollipop-tadpole}. The root of~$K_a$ is said to be the \emph{center} of $K_a^l$.
For any $0\le k\le a-1$, the \emph{melting lollipop} $K_a^l(k)$ is the graph obtained by removing $k$ edges from the clique $K_a$ that are incident with the center. Lollipops are the melting lollipops with $k=0$. \citet[Theorem 4.9]{HNY20} established the $e$-positivity of melting lollipops using generating function techniques. \citet{Tom25} obtained a positive $e_I$-expansion for melting lollipops.

\begin{theorem}[Melting lollipops, \citeauthor{Tom25}]\label{thm:lollipop.melting}
Let $n=a+l$, where $a\ge 2$ and $l\ge 0$. For any $0\le k\le a-1$,
\[
\frac{X_{K_a^l(k)}}{(a-2)!}
=
\sum_{I\vDash n,\
i_{-1}=a-1
}
k
w_{I\backslash i_{-1}}
e_I
+
\sum_{I\vDash n,\
i_{-1}\ge a
}
(a-k-1)
w_I
e_I.
\]
In particular, we have $X_{K_a^l}=(a-1)!\sum_{\substack{I\vDash n,\ i_{-1}\ge a}}w_I e_I$.
\end{theorem}

A \emph{KPK graph} is a $K$-chain of the form $P^k(K_a,K_b)$. It has order $a+b+k-1$. \citet[Theorem~3.6]{WZ25} obtained a positive $e_I$-expansion for the chromatic symmetric function of KPKs.

\begin{theorem}[KPKs, \citeauthor{WZ25}]\label{thm:KPK}
Let $n=a+b+l-1$, where $a,b\ge 1$ and $l\ge 0$. Then
\[
\frac{X_{P^l(K_a,K_b)}}
{(a-1)!\
(b-1)!}
=
\sum_{
I\vDash n,\
i_{-1}\ge a,\
i_1\ge b}
w_I
e_I
+
\sum_{
I\vDash n,\
i_{-1}\ge a,\
i_1\le b-1< i_2}
(i_2-i_1)
\prod_{j\ge 3}(i_j-1)
e_I.
\]
\end{theorem}

We simplify the specialization $b=3$ here. It will be used in \cref{sec:path.tw,sec:lollipop.tw}.
\begin{corollary}\label{cor:KPK3=KPC3}
Let $n=a+l+2$, where $a\ge 3$ and $l\ge 0$. Then
\[
\frac{X_{P^l(K_a,\,C_3)}}
{2(a-1)!}
=
ne_n
+(n-2)e_{1(n-1)}
+(n-4)e_{2(n-2)}
+
\sum_{\substack{
I\vDash n,\
a\le i_{-1}\le n-3,\
i_2\ge 3
}}
w_Ie_I.
\]
\end{corollary}
\begin{proof}
Let $t_I=\prod_{j\ge 3}(i_j-1)$. Taking $b=3$ in \cref{thm:KPK}, we obtain
\begin{multline*}
\frac{X_{P^{l}(K_a,K_3)}}{2(a-1)!}
=
\sum_{
I\vDash n,\
i_{-1}\ge a,\
i_1\ge 3}
w_Ie_I
+
\sum_{
I\vDash n,\
i_{-1}\ge a,\
i_1\le 2< i_2}
(i_2-i_1)
t_I
e_I\\
=
\sum_{
\substack{
I\vDash n\\
i_1\ge 3,\
i_{-1}\ge a
}}
w_I
e_I
+
\sum_{
\substack{
I\vDash n,\
\ell(I)\ge 2\\
i_1=1,\
i_2\ge 3,\
i_{-1}\ge a
}}
w_I
e_I
+
\sum_{
\substack{
I\vDash n,\
\ell(I)\ge 3\\
i_1=2,\
i_2\ge 3,\
i_{-1}\ge a
}}
(i_2-2)
t_I
e_I
+
(n-4)
e_{2(n-2)}.
\end{multline*}
Thus the desired formula is equivalent to
\[
\sum_{
\substack{
I\vDash n\\
i_1\ge 3,\
i_{-1}\ge a
}}
w_I
e_I
+
\sum_{
\substack{
I\vDash n,\
\ell(I)\ge 2\\
i_1=1,\
i_2\ge 3,\
i_{-1}\ge a
}}
w_I
e_I
+
\sum_{
\substack{
I\vDash n,\
\ell(I)\ge 3\\
i_1=2,\
i_2\ge 3,\
i_{-1}\ge a
}}
(i_2-2)
t_I
e_I
=
\sum_{\substack{
I\vDash n,\
i_{-1}\ge a,\
i_{-1}\ne n-2\\
I=n \text{ or }i_2\ge 3
}}
w_I
e_I.
\]
Let $\mathcal A=\{I\vDash n\colon \ell(I)\ge 2,\ i_{-1}\ge a\}$. Then the desired identity above can be recast as
\[
\sum_{
I\in\mathcal A,\
i_1\ge 3
}
w_I
e_I
+
\sum_{
I\in\mathcal A,\
\ell(I)\ge 3,\
i_1=2,\
i_2\ge 3
}
(i_2-2)
t_I
e_I
=
\sum_{
I\in\mathcal A,\
i_1\ge 2,\
i_2\ge 3,\
i_{-1}\ne n-2
}
w_I
e_I.
\]
Since $a\ge 3$, both sides with the restriction $\ell(I)=2$ have the same value $\sum_{\substack{I=i_1i_2\in\mathcal A,\ i_1\ge 3}}w_Ie_I$. Let $\mathcal B=\{I\vDash n\colon \ell(I)\ge 3,\ i_{-1}\ge a\}$. After restricting to $\ell(I)\ge3$, the desired identity becomes
\[
\sum_{
I\in\mathcal B,\
i_1\ge 3
}
w_I
e_I
+
\sum_{
I\in\mathcal B,\
i_1=2,\
i_2\ge 3
}
(i_2-2)
t_I
e_I
=
\sum_{
I\in\mathcal B,\
i_1=2,\
i_2\ge 3
}
w_I
e_I
+
\sum_{
I\in\mathcal B,\
i_1\ge 3,\
i_2\ge 3
}
w_I
e_I,
\]
or equivalently,
\[
\sum_{
I\in\mathcal B,\
i_1\ge 3,\
i_2=2
}
w_I
e_I
=
\sum_{
I\in\mathcal B,\
i_1=2,\
i_2\ge 3
}
i_2
t_I
e_I.
\]
This identity follows by exchanging $i_1$ and $i_2$ in every composition indexed on either side. This completes the proof.
\end{proof}

A \emph{PKP graph} is a $K$-chain of the form $P_{g+1}+K_a+P_{h+1}$. It has order $g+a+h$. A \emph{KKP graph} is a $K$-chain of the form $K_a+K_b+P_{h+1}$. It has order $a+b+h-1$. \citet[Theorems 5.3 and~5.4]{QTW26} obtained positive $e_I$-expansions for the chromatic symmetric functions of these graphs. Let
\begin{align*}
f_1(I,b)
&=(b-1)w_I e_I,\\
f_2(I,b)
&=(b-2)
i_{-1}
w_{I\backslash i_{-1}}
e_I,
\quad\text{and}\\
f_3(I,b)
&=(i_{-1}-b+1)
w_{I\backslash i_{-1}}
e_I.
\end{align*}

For any composition $I=i_1i_2\dotsm\vDash n$ and any integer $0\le a\le n$, define
\begin{align}
\label{def:sigma-Theta.+}
\sigma_I(a)
&=\min\{i_1+\dots+i_k\colon
0\le k\le \ell(I),\
i_1+\dots+i_k\ge a\},
\quad
\Theta_I(a)
=\sigma_I(a)-a,\quad\text{and}\\
\label{def:sigma-Theta.-}
\sigma_I^-(a)
&=\max\{i_1+\dots+i_k\colon
0\le k\le \ell(I),\
i_1+\dots+i_k\le a\},
\quad
\Theta_I^-(a)
=a-\sigma_I^-(a).
\end{align}
It is straightforward to see that $\Theta_I(a)\ge 0$ and $\Theta_I^-(a)\ge 0$. Moreover, these functions satisfy
\[
\Theta_I^-(a)
=\Theta_{\overline{I}}(n-a)
\quad\text{and}\quad
\sigma_K^-(k_1+a)
=\sigma_{K\backslash k_1}^-(a)
+k_1.
\]

\begin{theorem}[PKPs and KKPs, \citeauthor{QTW26}]\label{thm:KPK:PKP}
We have the following.
\begin{enumerate}
\item
Let $n=g+h+a$, where $g,h\ge 0$ and $a\ge 2$. Then
\[
\frac{X_{P_{g+1}+K_a+P_{h+1}}}
{(a-2)!}
=
(a-1)e_n
+
\sum_{
I\vDash n,\
\Theta_I(h+1)\ge a-1}
f_2(I,a)
+
\sum_{
I\vDash n,\
i_{-1}\ge a-1}
f_3(I,a).
\]
\item
Let $n=a+b+h-1$, where $a\ge 1$, $b\ge 2$ and $h\ge 0$. Then
\[
\frac{X_{K_a+K_b+P_{h+1}}}
{(a-1)!(b-2)!}
=
\sum_{\substack{
I\vDash n\\
i_{-1}\ge n-h}}
f_1(I,b)
-
\sum_{\substack{
I\vDash n,\
i_{-1}+i_{-2}\ge n-h\\
i_{-1}\le \min(a-1,\,b-2)}}
f_3(I,b)
+
\sum_{
\substack{
I\vDash n,\
i_{-1}+i_{-2}\ge n-h\\
\max(a,\,b)\le
i_{-1}\le n-h-1}}
f_3(I,b).
\]
\end{enumerate}
\end{theorem}
In fact, \citet[Theorem 5.3]{QTW26} proved
\[
\frac{X_{P_{g+1}+K_a+P_{h+1}}}
{(a-2)!}
=
\sum_{
\substack{
I\vDash n,\ 
i_{-1}\ge a-1\\
\Theta_I(h+1)\ge a-1}}
f_1(I,a)
+
\sum_{
\substack{
I\vDash n,\
i_{-1}\le a-2\\
\Theta_I(h+1)\ge a-1}}
f_2(I,a)
+
\sum_{
\substack{I\vDash n,\
i_{-1}\ge a-1\\
\Theta_I(h+1)\le a-2}}
f_3(I,a).
\]
We observe that
\begin{equation}\label{f123}
f_1(I,a)
-f_2(I,a)
-f_3(I,a)
=\begin{dcases*}
(a-1)e_n,
& if $I=n$,\\
0,
& otherwise.
\end{dcases*}
\end{equation}
This identity allows us to recast their formula as the form in \cref{thm:KPK:PKP}.

\citet[Theorem 3.3]{QTW26} also treated an analogous family of graphs called the \emph{KPC graphs}, which are of the form $P^l(K_a,C_c)$. When $a=1$, KPC graphs reduce to \emph{tadpoles}; we adopt the notation $C_c^l=P^l(K_1,C_c)$ for compactness, see the right illustration in \cref{fig:lollipop-tadpole}. Here the root of~$C_c$ is said to be the \emph{center} of $C_c^l$.

\begin{theorem}[KPCs, \citeauthor{QTW26}]\label{thm:KPC}
Let $n=a+l+c-1$, where $a\ge 1$, $l\ge 0$ and $c\ge 2$. Then
\[
X_{P^l(K_a,\,C_c)}
=(a-1)!\sum_{I\vDash n}
c_I
w_I
e_I,
\]
where
\[
c_I=\begin{cases*}
0, 
& if $\ell(I)\ge 2$ and $i_2<a$,\\
\displaystyle
i_2-a-l
+\frac{i_2-i_1}{i_2-1},
& if $\ell(I)\ge 2$, $i_1\le a-1$, and $i_2\ge a+l$,\\
\Theta_I(a+l),
& otherwise.
\end{cases*}
\]
In particular, we have $X_{C_c^l}=\sum_{I\vDash n}\Theta_I(l+1)w_Ie_I$.
\end{theorem}

The formula for tadpoles $C_c^l$ was first known by \citet[Theorem~3.2]{WZ25}. Using \cref{thm:KPC}, one may show \cref{cor:KPK3=KPC3} another way as below.

\begin{proof}[Another proof of \cref{cor:KPK3=KPC3}]
Let $G=P^l(K_a,\,C_3)$. Taking $c=3$ in \cref{thm:KPC}, we find $c_n=2$. Suppose that $i_1\le a-1$, $i_2\ge n-2$ and $w_I>0$. Then $I=2(n-2)$ or $I=1(n-1)$. One may compute that $c_{2(n-2)}=(n-4)/(n-3)$ and $c_{1(n-1)}=2$. In the remaining case, we have $i_2\ge a$ and either $i_1\ge a$ or $i_2\le n-3$. Since $a\ge 3$, this condition is equivalent to saying that $a\le i_2\le n-3$. Therefore,
\[
\frac{X_G}
{(a-1)!}
=2ne_n
+2(n-4)e_{2(n-2)}
+2(n-2)e_{1(n-1)}
+\sum_{I\vDash n,\
a\le i_2\le n-3}
\Theta_I(n-2)
w_I
e_I.
\]
If $\Theta_I(n-2)w_I\ne 0$, then $i_{-1}\ge 3$ and $\Theta_I(n-2)=2$. It follows that
\[
\frac{X_G}
{2(a-1)!}
=ne_n
+(n-4)e_{2(n-2)}
+(n-2)e_{1(n-1)}
+\sum_{I\vDash n,\
a\le i_2\le n-3,\
i_{-1}\ge 3}
w_I
e_I.
\]
In the last sum, one may exchange $i_2$ and $i_{-1}$ without loss of generality. The terms indexed by $I=n$ and $I=1(n-1)$ have the same form $w_Ie_I$. This completes the proof.
\end{proof}


\section{A positive \texorpdfstring{$e_I$}{eI}-expansion for KPKP graphs}\label{sec:KPKP}

A \emph{KPKP graph} is a graph of the form $K_a+P_{g+1}+K_b+P_{h+1}$, denoted $P^g(K_a,\,K_b^h)$; see \cref{fig:KPKP}. The two attachments to the middle clique $K_b$ use distinct vertices: the left path meets one vertex of $K_b$, while the right path meets the center of the lollipop $K_b^h$. The graph has order $a+g+b+h-1$ and is a $K$-chain. This section derives a positive $e_I$-expansion for its chromatic symmetric function.

\begin{theorem}[KPKP graphs]\label{thm:KPKP}
Let $n=a+g+b+h-1$, where $g,h\ge 0$, $a\ge 1$ and $b\ge 2$. Then we have 
the positive $e_I$-expansion
\begin{multline*}
\frac{X_{P^g(K_a,\,K_b^h)}}
{(a-1)!(b-2)!}
=
(b-1)ne_n
+
\sum_{\substack{
K\vDash n,\
\Theta_K(h+1)\ge b-1\\
k_{-1}+k_{-2}\le n-h-1\\
k_{-1}\ge b-1,\
k_{-2}\ge a
}}
f_1
+
\sum_{\substack{
K\vDash n,\
\Theta_K(h+1)\ge b-1\\
k_{-1}+k_{-2}
\ge n-h\\
k_{-1}
\ge \max(a,\,b-1)
}}
f_1
\\
+
\sum_{
\substack{
K\vDash n,\
\Theta_K(h+1)\ge b-1\\
k_{-1}\le b-2,\
k_{-2}\ge a\\
k_{-1}\ge a
\text{ or }
k_{-1}+k_{-2}\le n-h-1
}}
f_2
-
\sum_{
\substack{
K\vDash n,\
\Theta_K(h+1)\ge b-1\\
k_{-1}+k_{-2}
\ge n-h\\
k_{-1}
\le 
\min(a-1,\,b-2)
}}
f_3
+
\sum_{\substack{
K\vDash n,\
\Theta_K(h+1)\le b-2\\
k_{-1}\ge b-1\\
k_{-1}+k_{-2}\ge n-h
\text{ or }
k_{-2}\ge a
}}
f_3,
\end{multline*}
where $w$ and $\Theta$ are defined by \cref{def:w,def:sigma-Theta.+} respectively, and
\[
f_1
=(b-1)
w_K
e_K,
\qquad
f_2
=(b-2)
k_{-1}
w_{K\backslash k_{-1}}
e_K,
\qquad\text{and}\quad
f_3
=(k_{-1}-b+1)
w_{K\backslash k_{-1}}
e_K.
\]
\end{theorem}

\begin{proof}
The expansion is a positive one because (i) the fourth sum has only nonnagative coefficients since $f_3\le0$ by $k_{-1}\le b-2$, and (ii) the final sum has only nonnegative coefficients since $f_3\ge0$ by $k_{-1}\ge b-1$. 

Set $G=P^g(K_a,\,K_b^h)$. Taking $H=K_b^h$ in \cref{fml:KPG}, we obtain
\[
\frac{X_G}
{(a-1)!}
=
\sum_{l=0}^{a-1}
(1-l)e_l
X_{K_b^{a+g-1-l,\,h}}.
\]
In view of the formula for PKPs in \cref{thm:KPK:PKP}, we define
\begin{align*}
\mathcal S_{n1}
&=\{K\vDash n\colon 
\Theta_K(h+1)\ge b-1,\
k_{-1}\ge b-1\},\\
\mathcal S_{n2}
&=\{K\vDash n\colon 
\Theta_K(h+1)\ge b-1,\
k_{-1}\le b-2\},\\
\mathcal S_{n3}
&=\{K\vDash n\colon 
\Theta_K(h+1)\le b-2,\
k_{-1}\ge b-1\}.
\end{align*}
Then $X_G=(a-1)!(b-2)!(X_1+X_2+X_3)$, where
\begin{align*}
X_i
&=
\sum_{K\in \mathcal S_{ni}}
f_i(K,b)
-
\sum_{l=1}^{a-1}
\sum_{I\in\mathcal S_{(n-l)i}}
(l-1)e_l
f_i(I,b).
\end{align*}
First, we express the double sum in each $X_i$ by a single sum. For any $I\in\mathcal S_{(n-l)i}$, we consider the composition $K=Il$. Then
\[
(l-1)e_l
f_1(I,b)
=(l-1)e_l(b-1)w_Ie_I
=(b-1)w_Ke_K
=f_1(K,b).
\]
For brevity, we write $f_i=f_i(K,b)$ for $i=1,2,3$ when $K\vDash n$ is clear from context. Then
\begin{equation}\label{pf:KPKP.A1}
X_1
=
\sum_{K\in\mathcal S_{n1}}
f_1
-
\sum_{\substack{
K\vDash n,\
\Theta_K(h+1)\ge b-1\\
k_{-2}\ge b-1,\
k_{-1}\le a-1
}}
f_1.
\end{equation}
We proceed according to the number $j(K)$ such that $\sigma_K(h+1)=k_1+\dots+k_{-j(K)}$. 
\begin{enumerate}
\item
If $j(K)\ge 3$, i.e., if $k_{-1}+k_{-2}\le n-h-1$, then we can exchange $k_{-1}$ and $k_{-2}$ in the negative sum in \cref{pf:KPKP.A1}. Thus the difference for compositions $K$ with $j(K)\ge 3$ is
\[
\sum_{\substack{
K\in\mathcal S_{n1}\\
j(K)\ge 3
}}
f_1
-
\sum_{\substack{
K\vDash n,\
\Theta_K(h+1)\ge b-1\\
k_{-2}\ge b-1,\
k_{-1}\le a-1\\
j(K)\ge 3
}}
f_1
=
\sum_{\substack{
K\in\mathcal S_{n1}\\
k_{-1}+k_{-2}\le n-h-1
}}
f_1
-
\sum_{\substack{
K\in\mathcal S_{n1},\
k_{-2}\le a-1\\
k_{-1}+k_{-2}\le n-h-1
}}
f_1
=
\sum_{\substack{
K\in\mathcal S_{n1},\
k_{-2}\ge a\\
k_{-1}+k_{-2}\le n-h-1
}}
f_1.
\]
\item
Suppose that $j(K)\le 2$. When $\ell(K)\ge 2$, this premise is equivalent to saying that $k_{-1}+k_{-2}\ge n-h$. For any composition $K$ that appears in the negative sum of \cref{pf:KPKP.A1}, we have $\ell(K)\ge 2$ and
\[
k_{-2}
=(k_{-1}+k_{-2})-k_{-1}
\ge 
(n-h)-(a-1)
=g+b
>b-1.
\]
In other words, the restriction $k_{-2}\ge b-1$ in \cref{pf:KPKP.A1} is redundant. Since any function $F(s)$ defined on integers satisfies
\[
\sum_{s\ge u}F(s)
-\sum_{s\le v}F(s)
=\sum_{s\ge \max(u,\,v+1)}F(s)
-\sum_{s\le \min(u-1,\,v)}F(s)
\]
for any integers $u$ and $v$, we deduce that the difference for compositions~$K$ with $j(K)\le 2$ is
\begin{align*}
\sum_{
\substack{
K\in\mathcal S_{n1}\\
j(K)\le 2
}}
f_1
-
\sum_{\substack{
K\vDash n,\
\Theta_K(h+1)\ge b-1\\
k_{-2}\ge b-1,\
k_{-1}\le a-1\\
j(K)\le 2
}}
f_1
&=
f_1(n,b)
+
\sum_{\substack{
K\in\mathcal S_{n1}\\
k_{-1}+k_{-2}
\ge n-h}}
f_1
-
\sum_{\substack{
K\vDash n,\
\Theta_K(h+1)\ge b-1\\
k_{-1}
\le a-1\\
k_{-1}+k_{-2}
\ge n-h
}}
f_1\\
&=
(b-1)ne_n
+
\sum_{\substack{
K\vDash n,\
\Theta_K(h+1)\ge b-1\\
k_{-1}+k_{-2}
\ge n-h\\
k_{-1}
\ge \max(a,\,b-1)
}}
f_1
-
\sum_{\substack{
K\vDash n,\
\Theta_K(h+1)\ge b-1\\
k_{-1}+k_{-2}
\ge n-h\\
k_{-1}
\le 
\min(a-1,\,b-2)
}}
f_1.
\end{align*}
\end{enumerate}

Next, we handle $X_2$ in the same spirit. For any $I\in\mathcal S_{(n-l)2}$, we have $\ell(I)\ge 2$ and
\[
(l-1)e_l
f_2(I,b)
=(l-1)
e_l
(b-2)
i_{-1}
w_{I\backslash i_{-1}}
e_I
=(b-2)k_{-2}w_{K\backslash k_{-2}}e_K.
\]
It follows that
\begin{equation}\label{pf:KPKP.A2}
X_2
=
\sum_{K\in\mathcal S_{n2}}
f_2
-
\sum_{
\substack{
K\vDash n,\
\Theta_K(h+1)\ge b-1\\
k_{-2}\le b-2,\
k_{-1}\le a-1
}}
(b-2)
k_{-2}
w_{K\backslash k_{-2}}
e_K.
\end{equation}
Let $K$ be a composition that appears in the negative sum of \cref{pf:KPKP.A2}. Then
\[
k_1+\dots+k_{-3}
=n-k_{-1}-k_{-2}
\ge n-(b-2)-(a-1)=g+h+2>h+1.
\]
Thus $j(K)\ge 3$, which implies that the $\Theta$-restriction is independent of $k_{-1}$ and $k_{-2}$. This allows us to exchange $k_{-1}$ and $k_{-2}$ in the negative sum. Therefore,
\[
X_2
=
\sum_{
\substack{
K\vDash n,\
\Theta_K(h+1)\ge b-1\\
k_{-1}\le b-2
}}
f_2
-
\sum_{
\substack{
K\vDash n,\
\Theta_K(h+1)\ge b-1\\
k_{-1}\le b-2,\
k_{-2}\le a-1
}}
f_2
=
\sum_{
\substack{
K\vDash n,\
\Theta_K(h+1)\ge b-1\\
k_{-1}\le b-2,\
k_{-2}\ge a
}}
f_2.
\]

Finally, we deal with $X_3$ in the same fashion. For any $I\in\mathcal S_{(n-l)3}$, we have $\ell(K)\ge 2$ and
\[
(l-1)e_l
f_3(I,b)
=(l-1)
e_l
(i_{-1}-b+1)
w_{I\backslash i_{-1}}
e_I
=(k_{-2}-b+1)
w_{K\backslash k_{-2}}e_K.
\]
It follows that
\begin{equation}\label{pf:KPKP.A3}
X_3
=
\sum_{
K\in\mathcal S_{n3}
}
f_3
-
\sum_{\substack{
K\vDash n,\
\Theta_K(h+1)\le b-2\\
k_{-2}\ge b-1,\
k_{-1}\le a-1
}}
(k_{-2}-b+1)
w_{K\backslash k_{-2}}
e_K.
\end{equation}
\begin{enumerate}
\item
For compositions $K$ with $j(K)\ge 3$, we can exchange $k_{-1}$ and $k_{-2}$ in the negative sum of \cref{pf:KPKP.A3}, and obtain
\[
\sum_{
\substack{
K\in\mathcal S_{n3}\\
j(K)\ge 3
}}
f_3
-
\sum_{\substack{
K\vDash n,\
\Theta_K(h+1)\le b-2\\
k_{-2}\ge b-1,\
k_{-1}\le a-1,\
j(K)\ge 3
}}
(k_{-2}-b+1)
w_{K\backslash k_{-2}}
e_K
=
\sum_{\substack{
K\in\mathcal S_{n3}, k_{-2}\ge a\\
k_{-1}+k_{-2}\le n-h-1
}}
f_3.
\]
\item
If $j(K)\le 2$, i.e., if $k_{-1}+k_{-2}\ge n-h$, then
\[
b-2
\ge\Theta_K(h+1)
\ge k_1+\dots+k_{-2}-(h+1)
=n-h-1-k_{-1},
\]
i.e., $k_{-1}\ge n-h-b+1=a+g>a-1$. Thus this case does not happen for the negative sum. The positive sum for these compositions $K$ reduces to
\[
\sum_{
\substack{
K\in\mathcal S_{n3}\\
j(K)\le 2
}}
f_3
-
\sum_{\substack{
K\vDash n,\
\Theta_K(h+1)\le b-2\\
k_{-2}\ge b-1,\
k_{-1}\le a-1,\
j(K)\le 2
}}
(k_{-2}-b+1)
w_{K\backslash k_{-2}}
e_K
=
\sum_{\substack{
K\in\mathcal S_{n3}\\
k_{-1}+k_{-2}\ge n-h
}}
f_3.
\]
\end{enumerate}
Adding up $X_1$, $X_2$ and $X_3$, we obtain
\begin{multline*}
\frac{X_G}
{(a-1)!(b-2)!}
=
\sum_{\substack{
K\in\mathcal S_{n1}\\
k_{-1}+k_{-2}\le n-h-1\\
k_{-2}\ge a
}}
f_1
+
(b-1)ne_n
+
\sum_{\substack{
K\vDash n,\
\Theta_K(h+1)\ge b-1\\
k_{-1}+k_{-2}
\ge n-h\\
k_{-1}
\ge \max(a,\,b-1)
}}
f_1
\\
-
\sum_{\substack{
K\vDash n,\
\Theta_K(h+1)\ge b-1\\
k_{-1}+k_{-2}
\ge n-h\\
k_{-1}
\le 
\min(a-1,\,b-2)
}}
f_1
+
\sum_{
\substack{
K\vDash n,\
\Theta_K(h+1)\ge b-1\\
k_{-1}\le b-2,\
k_{-2}\ge a
}}
f_2
+
\sum_{\substack{
K\in\mathcal S_{n3}\\
k_{-1}+k_{-2}\ge n-h
\text{ or }
k_{-2}\ge a
}}
f_3.
\end{multline*}
We observe that every composition $K$ that appears in the negative sum above appears in the $f_2$-sum, since
\[
k_{-2}
=(k_{-1}+k_{-2})-k_{-1}
\ge (n-h)-(b-2)
=a+g+1
>a.
\]
Therefore, by \cref{f123}, one may merge these two sums as
\[
\sum_{
\substack{
K\vDash n,\
\Theta_K(h+1)\ge b-1\\
k_{-1}\le b-2,\
k_{-2}\ge a
}}
f_2
-
\sum_{\substack{
K\vDash n,\
\Theta_K(h+1)\ge b-1\\
k_{-1}+k_{-2}
\ge n-h\\
k_{-1}
\le 
\min(a-1,\,b-2)
}}
f_1
=
\sum_{
\substack{
K\vDash n,\
\Theta_K(h+1)\ge b-1\\
k_{-1}\le b-2,\
k_{-2}\ge a\\
k_{-1}\ge a
\text{ or }
k_{-1}+k_{-2}\le n-h-1
}}
f_2
-
\sum_{\substack{
K\vDash n,\
\Theta_K(h+1)\ge b-1\\
k_{-1}+k_{-2}
\ge n-h\\
k_{-1}
\le 
\min(a-1,\,b-2)
}}
f_3.
\]
Substituting it into the expression for $X_G$, we obtain the desired positive $e_I$-expansion.
\end{proof}

\Cref{thm:KPKP} has the following particular cases.
\begin{enumerate}
\item
When $b=2$, the graph $P^g(K_a,\,K_2^h)$ reduces to the lollipop $K_a^{g+h+1}$, and the formula in \cref{thm:KPKP} reduces to the one in \cref{thm:lollipop.melting}.
\item
When $h=0$, the graph $P^g(K_a,\,K_b^0)$ reduces to the KPK graph $P^g(K_a,\,K_b)$, and the formula in \cref{thm:KPKP} can be reduced to the one in \cref{thm:KPK} by merging suitable terms.
\item
When $g=0$, the graph $P^0(K_a,\,K_b^h)$ reduces to a KKP graph, and the formula in \cref{thm:KPKP} can be reduced to the one in \cref{thm:KPK:PKP} by merging suitable terms.
\item
When $a=1$ or $a=2$, the graph $P^g(K_a,\,K_b^h)$ reduces to a PKP graph, and the formula in \cref{thm:KPKP} reduces to the one in \cref{thm:KPK:PKP}.
\end{enumerate}

The specialization $b=3$ will be used in the proof of \cref{thm:lollipop.tw}. We simplify it here.

\begin{corollary}\label{cor:KPKP.m=3}
Let $n=a+g+h+2$, where $a\ge 1$ and $g,h\ge 0$. Then
\[
\frac{X_{P^g(K_a,\,K_3^h)}}
{(a-1)!}
=
\sum_{\substack{
K\vDash n,\
\Theta_K(h+1)\ge 2\\
k_{-1}
\ge a}}
f_1
+
\sum_{
\substack{
K\vDash n,\
\Theta_K(h+1)\ge 2\\
k_{-1}=1,\
k_{-2}\ge a
}}
f_2
+
\sum_{\substack{
K\vDash n,\
\Theta_K(h+1)\le 1,\
k_{-1}\ge 2\\
k_{-1}+k_{-2}\ge n-h
\text{ or }
k_{-2}\ge a
}}
f_3,
\]
where $f_1=2w_Ke_K$, $f_2=k_{-1}w_{K\backslash k_{-1}}e_K$, and $f_3=(k_{-1}-2)w_{K\backslash k_{-1}}e_K$.
\end{corollary}
\begin{proof}
Let $G=P^g(K_a,\,K_3^h)$. When $a=1$, the desired formula coincides with the one for PKPs in \cref{thm:KPK:PKP}. Below we suppose that $a\ge 2$. Taking $b=3$ in \cref{thm:KPKP}, we obtain
\[
\frac{X_G}
{(a-1)!}
=A+B+
\sum_{\substack{
K\vDash n,\
\Theta_K(h+1)\le 1,\
k_{-1}\ge 2\\
k_{-1}+k_{-2}\ge n-h
\text{ or }
k_{-2}\ge a
}}
f_3,
\]
where
\begin{align*}
A&=
2ne_n
+
\sum_{\substack{
K\vDash n,\
\Theta_K(h+1)\ge 2\\
k_{-1}+k_{-2}\le n-h-1\\
k_{-1}\ge 2,\
k_{-2}\ge a
}}
f_1
+
\sum_{\substack{
K\vDash n,\
\Theta_K(h+1)\ge 2\\
k_{-1}+k_{-2}
\ge n-h\\
k_{-1}
\ge a
}}
f_1
=
\sum_{\substack{
K\vDash n,\
\Theta_K(h+1)\ge 2\\
k_{-1}
\ge a}}
f_1,
\quad\text{and}\\
B&=
\sum_{
\substack{
K\vDash n,\
\Theta_K(h+1)\ge 2\\
k_{-1}=1,\
a\le k_{-2}\le n-h-2
}}
f_2
-
\sum_{
\substack{
K\vDash n,\
\Theta_K(h+1)\ge 2\\
k_{-1}=1,\
k_{-2}
\ge n-h-1
}}
f_3
=
\sum_{
\substack{
K\vDash n,\
\Theta_K(h+1)\ge 2\\
k_{-1}=1,\
k_{-2}\ge a
}}
f_2.
\end{align*}
In fact, in the first sum in $A$, the condition $k_{-1}\ge 2$ can be removed, and the parts $k_{-1}$ and $k_{-2}$ are exchangeable; while in $B$, we have $f_2=-f_3$ when $k_{-1}=1$. This completes the proof.
\end{proof}

\section{Positive \texorpdfstring{$e_I$}{eI}-expansions for twinned paths, twinned cycles and twinned lollipops}\label{sec:twin}

For any rooted graph $(G,v)$, the \emph{twinned graph} of $G$ at $v$ is the graph obtained by adding a new vertex $v'$, and by connecting $v'$ with $v$ and with all neighbors of $v$. In this section, we give positive $e_I$-expansions for twinned paths, twinned cycles and twinned lollipops by using the composition method.

\subsection{Twinned paths}\label{sec:path.tw}
\begin{figure}[htbp]
\begin{tikzpicture}[scale=1.5, decoration=brace]
\node (2) at (1, 0) [ellipsis]{};
\node (2l) at (1-\eps, 0) [ellipsis]{};
\node (2r) at (1+\eps, 0) [ellipsis]{};
\node (6) at (5, 0) [ellipsis]{};
\node (6l) at (5-\eps, 0) [ellipsis]{};
\node (6r) at (5+\eps, 0) [ellipsis]{};
\draw[edge] (0, 0) -- ($(2l) - (\eps, 0)$);
\draw[edge] ($(2r) + (\eps, 0)$) -- (2, 0) -- (3, 0) -- (4, 0) -- (3, .6) -- (2, 0);
\draw[edge] (3, 0) -- (3, .6);
\draw[edge] (4, 0) -- ($(6l) - (\eps, 0)$);
\draw[edge] ($(6r) + (\eps, 0)$) -- (6, 0);
\node (1) at (0, 0) [ball]{};
\node (3) at (2, 0) [ball]{};
\node (4) at (3, 0) [ball]{};
\node (4u) at (3, .6) [ball]{};
\node (5) at (4, 0) [ball]{};
\node (7) at (6, 0) [ball]{};
\draw [decorate] (2,-.1) -- (0,-.1);
\node at (1,0)[below=7pt] {length $l$};
\draw [decorate] (6,-.1) -- (4,-.1);
\node at (5,0)[below=7pt] {length $n-l-1$};
\end{tikzpicture}
\caption{The twinned path $\mathrm{tw}_l(P_n)$.}
\label{fig:tw.path.interior}
\end{figure}
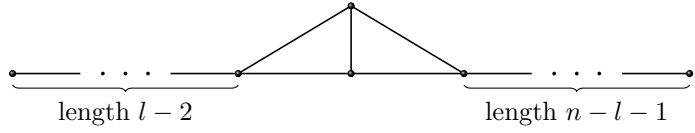

\citet{BCCCGKKLLS25} considered three ways of twinning a path, and showed that all the resulting graphs are $e$-positive by generating function techniques. We provide positive $e_I$-expansions for each of them. The first way is to twin a path at an end. This operation yields a lariat, for which one may find a positive $e_I$-expansion by taking $a=3$ and $l=n-3$ in the formula for lollipops in \cref{thm:lollipop.melting}. Precisely speaking, if we denote such a graph of order $n$ as $K_3^{n-3}$, then
\[
X_{K_3^{n-3}}
=2\sum_{I\vDash n,\ i_{-1}\ge 3}w_Ie_I,
\]
which is the specialization in \citet[Corollary~3.3]{WZ25}. The second way is to twin a path at both ends, which yields a barbell with two triangles. A positive $e_I$-expansion for this case is the specialization $a=3$ in \cref{cor:KPK3=KPC3}. The remaining case is to twin a path at an interior vertex. Let $G=\mathrm{tw}_l(P_n)$ be the graph that is obtained from the path $P_n=v_1\dotsm v_n$ by twinning at the vertex $v_{l+1}$, see \cref{fig:tw.path.interior}. It has order $n+1$.
\citet[Proposition~3.19]{BCCCGKKLLS25} showed that for $1\le l\le n-2$,
\begin{equation}\label{rec:tw.path}
X_G
=-2X_{P_l}X_{P_{n-l+1}}
+2e_1X_{P_n}
+4X_{P_{n+1}}
-2X_{P_{l+1}}X_{P_{n-l}}
+2e_2X_{P_l}X_{P_{n-l-1}}
-2X_{P_{l+2}}X_{P_{n-l-1}}.
\end{equation}

\begin{theorem}[Twinned path at an interior vertex]\label{thm:path.tw}
For $n\ge 3$ and $1\le l\le n-2$, the chromatic symmetric function of the twinned path $\mathrm{tw}_l(P_n)$ has the following positive $e_I$-expansion:
\begin{align*}
X_{\mathrm{tw}_l(P_n)}
&=2\sum_{I\vDash n,\ \Theta_I(l)\ge 3}w_Ie_{1I}
+2\sum_{I\in\mathcal W_n,\ \Theta_I(n-l-1)\ge 3}w_Ie_{1I}
+2\sum_{I\in\mathcal W_n}\brk3{1-\frac{2}{i_1}}w_Ie_{1I}\\
&\qquad
+2\sum_{K\in\mathcal W_{n+1},\ \Theta_K(l)\le 2}
\brk3{1-\frac{1}{\Theta_K(l+1+\Theta_K(l))}}w_Ke_K
+4\sum_{K\in\mathcal W_{n+1},\ \Theta_K(l)\ge 3}w_Ke_K,
\end{align*}
where $w$, $\mathcal W_n$ and $\Theta$ are defined by \cref{def:w,def:Cn,def:sigma-Theta.+} respectively.
\end{theorem}
\begin{proof}
Let $G=\mathrm{tw}_l(P_n)$ and let $h=l+1$. By \cref{rec:tw.path,prop:path},
\begin{equation}\label{pf:Pnl}
\frac{X_G}{2}
=2
\sum_{K\vDash n+1}
w_K 
e_K
+\sum_{K\vDash n}
w_K 
e_{1K}
+\sum_{
\substack{
I\vDash h-1\\
J\vDash n-h
}}
w_I
w_J
e_{I\!J2}
-\sum_{k=0}^2
\sum_{\substack{
I\vDash h-1+k\\
J\vDash n-h-k+2
}}
w_I
w_J
e_{I\!J}.
\end{equation}
It is clear that $[e_1^k]X_G=0$ for $k\ge 3$. Let
\[
\frac{X_G}{2}
=\sum_{k=0}^2
Y_k
e_1^k
\]
be the polynomial expansion in $e_1$. We proceed by computing each $Y_k$.

Extracting the coefficient of $e_1^2$ from \cref{pf:Pnl}, we obtain
\begin{align*}
Y_2
&=\sum_{K\vDash n-1}
w_{1K} 
e_K
+\sum_{\substack{
I\vDash h-2\\
J\vDash n-h-1
}}
w_{1I\!J}
e_{2I\!J}
-\sum_{k=0}^2
\sum_{\substack{
I\vDash h-2+k\\
J\vDash n-h-k+1
}}
w_{1I\!J}
e_{I\!J}\\
&=\sum_{K\vDash n-1}
w_{1K} 
e_K
+\sum_{
\substack{
K\vDash n-3\\
\Theta_K(h-2)=0}}
w_{1K}
e_{2K}
-\sum_{\substack{
K\vDash n-1\\
\Theta_K(h-2)\le 2}}
w_{1K}
e_K
-\sum_{\substack{
K=I2J\vDash n-1\\
I\vDash h-2}}
w_{1K}
e_{2K}\\
&=\sum_{K\vDash n-1,\
\Theta_K(h-2)\ge 3}
w_{1K} 
e_K.
\end{align*}
Second, extracting the coefficient of $e_1$ from \cref{pf:Pnl} (excluding the terms in $Y_2$), we obtain
\begin{align*}
Y_1
&=2
\sum_{K\in\mathcal W_n}
w_{1K} 
e_K
+\sum_{K\in\mathcal W_n}
w_K 
e_K\\
&\quad+
\sum_{
K=I\!J\in\mathcal W_{n-2},\
I\in\mathcal W_{h-1}}
w_K
e_{K2}
+\sum_{
K=J\!I\in\mathcal W_{n-2},\
I\in\mathcal W_{h-2}}
w_K
e_{K2}
-T_1-T_2,
\end{align*}
where
\begin{align*}
T_1
&=\sum_{k=0}^2
\sum_{
\substack{
K=I\!J\in\mathcal W_n\\
I\in \mathcal W_{h-1+k}}}
w_K
e_K
=\sum_{\substack{
K\in\mathcal W_n\\
\Theta_K(h-1)\le 2}}
w_K
e_K
+\sum_{\substack{
K=I2J\in\mathcal W_n\\
I\in\mathcal W_{h-1}}}
w_K
e_K,
\quad\text{and}\\
T_2
&=\sum_{k=0}^2
\sum_{\substack{
K=J\!I\in\mathcal W_n\\
I\in \mathcal W_{h-2+k}}}
w_K
e_K
=\sum_{\substack{
K\in\mathcal W_n\\
\Theta_K(n-h)\le 2}}
w_K
e_K
+\sum_{\substack{
K=I2J\in\mathcal W_n\\
I\in\mathcal W_{n-h}}}
w_K
e_K,
\end{align*}
Canceling the last sum in $T_1$ and the last sum in $T_2$ with the positive sums in the parenthesis of $Y_1$, we obtain
\[
Y_1
=2
\sum_{K\in\mathcal W_n}
w_{1K} 
e_K
+\sum_{K\in\mathcal W_n}
w_K 
e_K
-\sum_{
K\in\mathcal W_n,\
\Theta_K(h-1)\le 2}
w_K
e_K
-\sum_{
K\in\mathcal W_n,\
\Theta_K(n-h)\le 2}
w_K
e_K.
\]
Expressing the negative sums above by their ``$\Theta$-complements'', we infer that
\begin{align*}
Y_1
&=2
\sum_{K\in\mathcal W_n}
w_{1K} 
e_K
-\sum_{K\in\mathcal W_n}
w_K 
e_K
+\sum_{
K\in\mathcal W_n,\
\Theta_K(h-1)\ge 3}
w_K
e_K
+\sum_{
K\in\mathcal W_n,\
\Theta_K(n-h)\ge 3}
w_K
e_K\\
&=
\sum_{K\in\mathcal W_n}
\frac{k_1-2}{k_1}w_K 
e_K
+\sum_{
K\in\mathcal W_n,\
\Theta_K(h-1)\ge 3}
w_K
e_K
+\sum_{
K\in\mathcal W_n,\
\Theta_K(n-h)\ge 3}
w_K
e_K.
\end{align*}
Finally, extracting $Y_0$ from \cref{pf:Pnl} and treating the negative sums in the same way, we obtain
\begin{align*}
Y_0
&=2
\sum_{K\in\mathcal W_{n+1}}
w_K 
e_K
+\sum_{\substack{
K=I\!J\in\mathcal W_{n-1}\\
I\vDash h-1}}
w_I
w_J
e_{K2}
-\sum_{k=0}^2
\sum_{\substack{
K=I\!J\in\mathcal W_{n+1}\\
I\vDash h-1+k}}
w_I
w_J
e_K\\
&=2
\sum_{K\in\mathcal W_{n+1}}
w_K 
e_K
-
\sum_{
K=I\!J\in\mathcal W_{n+1},\
\Theta_K(h-1)\le 2,\
\abs{I}=\sigma_K(h-1)}
w_I
w_J
e_K\\
&=2
\sum_{\substack{
K\in\mathcal W_{n+1}\\
\Theta_K(h-1)\ge 3
}}
w_K 
e_K
+
\sum_{\substack{
K=I\!J\in\mathcal W_{n+1}\\
\Theta_K(h-1)\le 2,\
\abs{I}=\sigma_K(h-1)}}
(2w_K-w_Iw_J)
e_K,
\end{align*}
in which
\[
2w_K-w_Iw_J
=\brk3{1-\frac{1}{j_1-1}}
w_K.
\]
For any pair $(I,J)$ that appears in the last sum of $Y_0$,
\[
j_1
=\Theta_K\brk1{
h+\Theta_K(h-1)}
+1.
\]
Note that $\Theta_K(h+\Theta_K(h-1))\ge 1$, which implies $j_1\ge 2$. It follows that
\[
Y_0
=2
\sum_{
K\in\mathcal W_{n+1},\
\Theta_K(h-1)\ge 3
}
w_K
e_K
+\sum_{
K\in\mathcal W_{n+1},\
\Theta_K(h-1)\le 2}
\brk3{
1-\frac{1}{
\Theta_K(h+\Theta_K(h-1))}}
w_K
e_K.
\]
In conclusion,
\begin{multline*}
\frac{X_G}{2}
=\sum_{\substack{
K\vDash n-1\\
\Theta_K(h-2)\ge 3}}
w_{1K} 
e_{K11}
+
\sum_{K\in\mathcal W_n}
\frac{k_1-2}{k_1}w_K 
e_{1K}
+\sum_{\substack{
K\in\mathcal W_n\\
\Theta_K(h-1)\ge 3}}
w_K
e_{1K}
+\sum_{\substack{
K\in\mathcal W_n\\
\Theta_K(n-h)\ge 3}}
w_K
e_{1K}\\
+2
\sum_{\substack{
K\in\mathcal W_{n+1}\\
\Theta_K(h-1)\ge 3
}}
w_K
e_K
+\sum_{\substack{
K\in\mathcal W_{n+1}\\
\Theta_K(h-1)\le 2}}
\brk3{
1-\frac{1}{
\Theta_K(h+\Theta_K(h-1))}}
w_K
e_K,
\end{multline*}
in which the first sum and the third sum can be merged as desired.
\end{proof}

\subsection{Twinned cycles}\label{sec:cycle.tw}
\begin{figure}[htbp]
\centering
\begin{tikzpicture}
\draw (0,0) node {$C_n$};
\node[ball] (v) at (0: \r) {};
\node[ball] (u) at (30: \r) {};
\node[ball] (w) at (-30: \r) {};
\node[ball] (v') at (0: 1.5*\r) {};
\draw[edge] (-75: \r) arc (-75: 75: \r);
\draw[edge] (30: \r) -- (0: 1.5*\r);
\draw[edge] (0: \r) -- (0: 1.5*\r);
\draw[edge] (-30: \r) -- (0: 1.5*\r);
\end{tikzpicture}
\caption{The twinned cycle $\mathrm{tw}(C_n)$.}\label{fig:tw.cycle}
\end{figure}
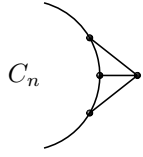

The \emph{twinned cycle} $G=\mathrm{tw}(C_n)$ is the graph of order $n+1$ obtained by twinning a vertex of the cycle $C_n$, see \cref{fig:tw.cycle}.
\citet[Lemma 3.26]{BCCCGKKLLS25} showed that for $n\ge 3$ and for the twinned cycle $G=\mathrm{tw}(C_n)$,
\begin{equation}\label{rec:Cnv}
X_G
=
4X_{C_{n+1}}
+2e_1X_{C_n}
-6X_{P_{n+1}}
+2e_2X_{P_{n-1}}.
\end{equation}
Based on \cref{rec:Cnv}, they confirmed the $e$-positivity of twinned cycles using generating function techniques. We now present a positive $e_I$-expansion for $X_{\mathrm{tw}(C_n)}$, which gives a direct proof of $e$-positivity.

\begin{theorem}[Twinned cycle]\label{thm:cycle.tw}
For any $n\ge 3$, we have
\[
X_{\mathrm{tw}(C_n)}
=
\smashoperator[r]{
\sum_{I\vDash n,\ i_1\ge 4
}}
2(i_1-3)
w_{1I}
e_{1I}
+
\smashoperator[r]{
\sum_{I\in\mathcal W_{n+1},\
i_1,i_{-1}\ge 3
}}
2(2i_1-5)
w_I
e_I
+\sum_{I\in\mathcal W_{n-1},\
i_1\ge 3}
4\brk3{i_1-3+\frac{1}{i_1}}
w_I
e_{I2},
\]
where $w$ and $\mathcal W_n$ are defined by \cref{def:w,def:Cn} respectively.
\end{theorem}
\begin{proof}
Let $G=\mathrm{tw}(C_n)$. By \cref{rec:Cnv,prop:cycle,prop:path},
\begin{equation}\label{pf:Cnv}
\frac{X_G}{2}
=2
\sum_{I\vDash n+1}
(i_1-1)
w_I
e_I
+\sum_{I\vDash n}
(i_1-1)
w_I
e_{I1}
-3\sum_{I\vDash n+1}
w_I
e_I
+\sum_{I\vDash n-1}
w_I
e_{I2}.
\end{equation}
Let $X_G/2=Y_0+Y_1e_1$ be the polynomial expansion of $X_G/2$ in $e_1$. We compute $Y_0$ and $Y_1$ separately.

First, extracting the $e_1$-coefficient from \cref{pf:Cnv}, we obtain
\[
Y_1
=\sum_{I\vDash n}
(i_1-1)
w_I
e_I
-3\sum_{I\vDash n}
w_{1I}
e_I
+\sum_{I\vDash n-2}
w_{1I}
e_{I2}.
\]
Since $(i_1-1)w_I=i_1w_{1I}$, we can simplify $Y_1$ as
\[
Y_1
=\sum_{I\vDash n}
(i_1-3)
w_{1I}
e_I
+\sum_{I\vDash n-2}
w_{1I}
e_{I2}
=\sum_{I\vDash n,\ i_1\ge 4}
(i_1-3)
w_{1I}
e_I.
\]
Second, extracting the terms of \cref{pf:Cnv} that do not contain $e_1$, we obtain
\begin{align*}
Y_0
&=
\sum_{I\in\mathcal W_{n+1}}
(2i_1-5)
w_I
e_I
+\sum_{I\in\mathcal W_{n-1}}
w_I
e_{I2}\\
&=
\sum_{I\in\mathcal W_{n+1},\
i_1\ge 3}
(2i_1-5)
w_I
e_I
+\sum_{I\in\mathcal W_{n-1}}
(w_I-w_{2I})
e_{I2},
\end{align*}
in which $w_I-w_{2I}=\frac{2-i_1}{i_1}w_I$. It follows that
\[
Y_0
=\sum_{I\in\mathcal W_{n+1},\
i_1,i_{-1}\ge 3}
(2i_1-5)
w_I
e_I
+\sum_{I\in\mathcal W_{n-1},\
i_1\ge 3}
(2i_1-5)
w_I
e_{I2}
+\sum_{I\in\mathcal W_{n-1}}
\frac{2-i_1}{i_1}
w_I
e_{I2}.
\]
In view of the factor $2-i_1$, we can regard $i_1\ge 3$ in the last sum. Therefore,
\[
Y_0
=\sum_{I\in\mathcal W_{n+1},\
i_1,i_{-1}\ge 3}
(2i_1-5)
w_I
e_I
+\sum_{I\in\mathcal W_{n-1},\
i_1\ge 3}
2\brk3{i_1-3+\frac{1}{i_1}}
w_I
e_{I2}.
\]
Combining $Y_1$ and $Y_0$, we obtain the desired formula.
\end{proof}

For $n=3$, $G=K_4$, only the second sum survives, and $X_{K_4}=2\cdotp 3\cdotp w_4e_4=24e_4$.
For example,
\begin{align*}
X_{\mathrm{tw}(C_4)}
&=
50e_5
+6e_{41}
+4e_{32},\\
X_{\mathrm{tw}(C_5)}
&=
84e_6
+16e_{51}
+20e_{42}
+12e_{33},
\quad\text{and}\\
X_{\mathrm{tw}(C_6)}
&=
126e_7
+30e_{61}
+44e_{52}
+66e_{43}
+6e_{421}
+4e_{322}.
\end{align*}

\subsection{Twinned lollipops}\label{sec:lollipop.tw}
Let $H$ be the graph obtained by twinning the lollipop $K_a^l$ at a vertex~$v^*$. Then $H$ has order $n=a+l+1$. There are three possible choices for $v^*$.
\begin{enumerate}
\item
If $v^*$ lies in the clique $K_a$, then $H=K_{a+1}^l$ is a lollipop or $H=K_{a+2}^{l-1}(a-1)$ is a melting lollipop; these cases are covered by \cref{thm:lollipop.melting}.
\item
If $v^*$ is the leaf, then $H=P^{l-1}(K_a,K_3)$ is a KPK graph, which is covered by \cref{cor:KPK3=KPC3}.
\item
Suppose that $v^*$ is an interior vertex in the path part. Let $h$ be the length of the path from $v^*$ to the leaf, see \cref{fig:lollipop.tw}. We write $H=\mathrm{tw}_h(K_a^l)$ in this case.
\end{enumerate}

\begin{theorem}[Twinned lollipop at an interior vertex of the path part]\label{thm:lollipop.tw}
Let $n=a+l+1$, where $a\ge 1$ and $l\ge 2$. For any $1\le h\le l-1$,
\begin{align*}
\frac{X_{\mathrm{tw}_h(K_a^l)}}
{2(a-1)!}
&=
\smashoperator[r]{
\sum_{\substack{
K\vDash n,\
k_{-1}
\ge a\\
\Theta_K(h)\ge 3
}}}
2
w_K
e_K
+
\smashoperator[r]{
\sum_{
\substack{
I\vDash n-1,\
i_{-1}\ge a\\
\Theta_I(h)\ge 3
}}}
w_I
e_{1I}
+
\smashoperator[r]{
\sum_{\substack{
K\vDash n,\
\Theta_K(h)\le 1,\
k_{-1}\ge 3\\
k_{-1}+k_{-2}\ge n-h+1
\text{ or }
k_{-2}\ge a
}}}
(k_{-1}-2)
w_{K\backslash k_{-1}}
e_K\\
&\qquad+
\sum_{
K\vDash n,\
k_{-1}\ge a,\
\Theta_K(h)=2,\
\Theta_K(h+3)\ge2
}
\frac{\Theta_K(h+3)-1}
{\Theta_K(h+3)}
w_K
e_K,
\end{align*}
where the functions $w$ and $\Theta$ are defined by \cref{def:w,def:sigma-Theta.+} respectively. As a consequence, all twinned lollipops are $e$-positive.
\end{theorem}

\begin{proof}
Let $G=\mathrm{tw}_h(K_a^l)$ and $g=l-h-1$. Applying \cref{prop:3del}, we obtain
\begin{equation}\label{pf:rec:TwinnedLollipop}
X_G
=2X_{P^{g+1}(K_a,\,K_3^{h-1})}
-
X_{K_3^{h-1}}
X_{K_a^g}.
\end{equation}
Write $f_i=f_i(K,3)$ for short. By \cref{cor:KPKP.m=3},
\[
\frac{X_{P^{g+1}(K_a,\,K_3^{h-1})}}
{(a-1)!}
=
\sum_{\substack{
K\vDash n,\
\Theta_K(h)\ge 2\\
k_{-1}
\ge a}}
f_1
+
\sum_{
\substack{
K\vDash n,\
\Theta_K(h)\ge 2\\
k_{-1}=1,\
k_{-2}\ge a
}}
f_2
+
\sum_{\substack{
K\vDash n,\
\Theta_K(h)\le 1,\
k_{-1}\ge 2\\
k_{-1}+k_{-2}\ge n-h+1
\text{ or }
k_{-2}\ge a
}}
f_3,
\]
in which
\[
\sum_{
\substack{
K\vDash n,\
\Theta_K(h)\ge 2\\
k_{-1}=1,\
k_{-2}\ge a
}}
f_2
=
\sum_{
\substack{
K\vDash n,\
\Theta_K(h)\ge 2\\
k_{-1}=1,\
k_{-2}\ge a
}}
w_{K\backslash k_{-1}}
e_K
=
\sum_{
\substack{
I\vDash n-1,\
\Theta_I(h)\ge 2\\
i_{-1}\ge a
}}
w_I
e_{1I}.
\]
For $a\ge2$, the following product formula follows from the lollipop specialization in \cref{thm:lollipop.melting}. For $a=1$, it follows from \cref{prop:path}, since $K_1^g=P_{g+1}$. Thus
\[
\frac{X_{K_3^{h-1}}X_{K_a^g}}
{2(a-1)!}
=
\smashoperator[r]{
\sum_{
\substack{
I\vDash h+2,\
i_{-1}\ge 3\\
J\vDash n-h-2,\
j_{-1}\ge a
}}}
w_I
w_J
e_{I\!J}
=
\smashoperator[r]{
\sum_{
\substack{
I\vDash n-1,\
\Theta_I(h)=2\\
i_{-1}\ge a
}}}
w_I
e_{1I}
+
\sum_{
\substack{
K\vDash n,\
\Theta_K(h)=2\\
k_{-1}\ge a,\
\Theta_K(h+3)\ge 1
}}
\frac{\Theta_K(h+3)+1}{\Theta_K(h+3)}
w_K
e_K,
\]
since $j_1=\Theta_K(h+3)+1$.
Substituting these two results into \cref{pf:rec:TwinnedLollipop}, we obtain
\[
\frac{X_G}
{2(a-1)!}
=A+B+\sum_{\substack{
K\vDash n,\
\Theta_K(h)\le 1,\
k_{-1}\ge 2\\
k_{-1}+k_{-2}\ge n-h+1
\text{ or }
k_{-2}\ge a
}}
f_3,
\]
where
\begin{align*}
A
&=
\sum_{\substack{
K\vDash n,\
\Theta_K(h)\ge 2\\
k_{-1}
\ge a}}
f_1
-
\sum_{
\substack{
K\vDash n,\
\Theta_K(h)=2\\
k_{-1}\ge a,\
\Theta_K(h+3)\ge 1
}}
\frac{\Theta_K(h+3)+1}{\Theta_K(h+3)}
w_K
e_K,\quad\text{and}\\
B
&=
\sum_{
\substack{
I\vDash n-1,\
\Theta_I(h)\ge 2\\
i_{-1}\ge a
}}
w_I
e_{1I}
-
\sum_{
\substack{
I\vDash n-1,\
\Theta_I(h)=2\\
i_{-1}\ge a
}}
w_I
e_{1I}
=
\sum_{
\substack{
I\vDash n-1,\
\Theta_I(h)\ge 3\\
i_{-1}\ge a
}}
w_I
e_{1I}.
\end{align*}
Recall that $f_1=2w_Ke_K$ and $n=a+g+h+2\ge h+3$. When $\Theta_K(h)=2$ and $f_1\ne0$, we have $\Theta_K(h+3)\ge 1$. Therefore,
\[
A
=
\sum_{\substack{
K\vDash n,\
\Theta_K(h)\ge 3\\
k_{-1}\ge a
}}
2w_Ke_K
+
\sum_{
\substack{
K\vDash n,\
\Theta_K(h)=2\\
k_{-1}\ge a,\
\Theta_K(h+3)\ge 1
}}
\frac{\Theta_K(h+3)-1}{\Theta_K(h+3)}
w_K
e_K,
\]
in which the restriction $\Theta_K(h+3)\ge 1$ can be improved to $\Theta_K(h+3)\ge 2$ without loss of generality. Substituting the functions $f_i$ by their definitions, we obtain the desired formula. The $e$-positivity of all twinned lollipops then follows from the formula and the arguments at the beginning of \cref{sec:lollipop.tw}.
\end{proof}

We remark that the formula in \cref{thm:lollipop.tw} does not hold for $h=0$, which can be seen by calculating the $e_n$-coefficient. In fact, the condition $h\ge 1$ is used when we invoke \cref{cor:KPKP.m=3,thm:lollipop.melting} in the proof above.

\begin{corollary}[Twinned path at an interior vertex]\label{cor:lollipop.tw}
Let $n\ge 3$. For any $1\le l\le n-2$,
\begin{align*}
X_{\mathrm{tw}_l(P_n)}
&=2\sum_{\substack{K\vDash n+1,\ k_{-1}\ge 3\\ \Theta_K(l)\le 1}}
(k_{-1}-2)w_{K\backslash k_{-1}}e_K
+2\sum_{\substack{K\vDash n+1,\ \Theta_K(l+3)\ge 2\\ \Theta_K(l)=2}}
\frac{\Theta_K(l+3)-1}{\Theta_K(l+3)}w_Ke_K\\
&\qquad
+4\sum_{K\vDash n+1,\ \Theta_K(l)\ge 3}w_Ke_K
+2\sum_{I\vDash n,\ \Theta_I(l)\ge 3}w_Ie_{1I},
\end{align*}
where the functions $w$ and $\Theta$ are defined by \cref{def:w,def:sigma-Theta.+} respectively.
\end{corollary}
\begin{proof}
Replacing the triple $(a,h,l)$ in \cref{thm:lollipop.tw} with $(1,l,n-1)$, and reflecting the path if necessary, the graph becomes $\mathrm{tw}_l(P_n)$, and the formula in \cref{thm:lollipop.tw} reduces to the desired one.
\end{proof}

The nonuniqueness of positive $e_I$-expansions is visible in $X_{\mathrm{tw}_1(P_6)}$: \cref{thm:path.tw} contains the contribution $12e_{142}$, whereas \cref{cor:lollipop.tw} distributes the same elementary-basis contribution as $4e_{124}+8e_{142}$.

\begin{figure}[htbp]
\begin{tikzpicture}
\node (1) at (0, 0) [ball]{};
\node (2) at (1, 0) [ball]{};
\node (3) at (2, 0) [ball]{};
\node (4) at (1, -.6) [ball]{};
\node (5) at (1, .6) [ball]{};
\node (6) at (3, 0) [ball]{};
\draw[edge] (0, 0) -- (1, 0) -- (2, 0) -- (1, .6) -- (0, 0);
\draw[edge] (1, 0) -- (1, .6);
\draw[edge] (0, 0) -- (1, -.6)  -- (2, 0) -- (3, 0);

\begin{scope}[xshift=5cm]
\node (1) at (0, 0) [ball]{};
\node (2) at (1, 0) [ball]{};
\node (3) at (2, 0) [ball]{};
\node (4) at (1, -.6) [ball]{};
\node (5) at (1, .6) [ball]{};
\node (6) at (3, 0) [ball]{};
\draw[edge] (0, 0) -- (1, -.6) -- (1, 0) -- (1, .6) -- (0, 0);
\draw[edge] (0, 0) -- (1, 0);
\draw[edge] (1, -.6) -- (2, 0)  -- (1, .6);
\draw[edge] (2, 0) -- (3, 0);
\end{scope}
\end{tikzpicture}
\caption{The twinned tadpoles $C_{4,1}^1$ and $C_{4,2}^1$.}
\label{fig:C41.tw}
\end{figure}
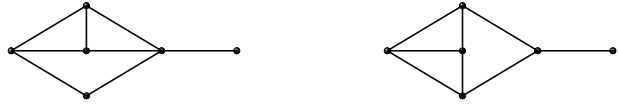

We remark that the twinning operation does not preserve the $e$-positivity if the lollipop is replaced with a tadpole. For instance, let $C_{4,i}^1$ be the graph obtained from the tadpole $C_4^1$ by twinning a vertex on $C_4$ that has distance $i$ with the center, see \cref{fig:C41.tw}.
Then $C_{4,1}^1$ is not $e$-positive while $C_{4,2}^1$ is $e$-positive, as
\begin{align*}
X_{C_{4,1}^1}
&=60e_6
+50e_{51} 
-4e_{42} 
+6e_{411} 
+6e_{33} 
+2e_{321}
\quad\text{and}\\
X_{C_{4,2}^1}
&=60e_6
+40e_{51} 
+12e_{42} 
+6e_{33} 
+2e_{321}.
\end{align*}

\section{Kayak paddles and infinity graphs}\label{sec:kayak}

A \emph{kayak paddle} is a graph of the form $P^l(C_a,C_b)$, where $a,b\ge 2$ and $l\ge 0$, see \cref{fig:kayak}. It has order $a+b+l-1$. The main result of this section is \cref{thm:kayak}. Its proof uses cyclic averaging and combinatorial injections.

\begin{theorem}[Kayak paddles]\label{thm:kayak}
Let $a,b\ge 3$, $l\ge 0$, and $n=a+b+l-1$. Then the chromatic symmetric function of the kayak paddle $P^l(C_a,C_b)$ has the following positive $e_I$-expansion:
\begin{align*}
X_{P^l(C_a,C_b)}
&=
\sum_{
\substack{
K\vDash n\\
k_1=1
}}
\Theta_K^-(a)g(K)
+
\sum_{
\substack{
K\vDash n,\
\Theta_K(a)\le l\\
2\le k_1\le l+1
}}
\Theta_K^-(k_1+a-1)
g(K)
\\
&\qquad
+
\sum_{
\substack{
K\vDash n,\
\Theta_K(a)\le l\\
l+2\le k_1\le l+a-1
}}
\brk1{
\Theta_K^-(a+l)+k_1-l-1}
g(K)
+
\sum_{
\substack{
K\vDash n\\
k_1\ge a+l
}}
(a-1)
g(K)
\\
&\qquad
+
\sum_{
\substack{
(I,J)\in \mathcal H\\
i_1\le l+1
\textnormal{ or }
i_1=j_1\\
\textnormal{or }
(i_1\ge l+2
\textnormal{ and }
\abs{I}\ge a+i_1-j_1)
}}
\brk3{
a-1-\abs{I}
+
\frac{j_1-i_1}
{j_1-1}}
g(I\!J)
\\
&\qquad
+
\sum_{
(I,J)\in\mathcal H,\
i_1>j_1
}
\brk4{
(i_1-j_1)
+
\brk1{
a-1-\abs{I}}
\brk3{
1+
\frac{j_1}{i_1}
\cdotp
\frac{i_1-1}{j_1-1}
}}
g(I\!J),
\end{align*}
where $g(K)=\Theta_K(a+l)w_Ke_K$ and $\mathcal H=\{(I,J)\colon I\!J\in\mathcal W_n,\ I\ne\emptyset,\ J\ne\emptyset,\ a+l-j_1+1\le \abs{I}\le a-1\}$. Here $w_K$, $\Theta_K$ and $\Theta_K^-$ are defined by \cref{def:w,def:sigma-Theta.+,def:sigma-Theta.-}, respectively.
\end{theorem}
\begin{proof}
Let $G=P^l(C_a,\,C_b)$. Taking $H=C_b$ in \cref{fml:CPG}, we obtain
\[
X_G
=(a-1)X_{C_b^{a+l-1}}
-\sum_{k=1}^{a-2}
X_{C_{a-k}}
X_{C_b^{k+l-1}}.
\]
Applying \cref{prop:cycle,thm:KPC}, we obtain
\[
X_G=(a-1)
\sum_{K\vDash n}
\Theta_K(a+l)
w_K
e_K
-\sum_{k=1}^{a-2}
\sum_{
I\vDash a-k,\
J\vDash n-\abs{I}}
(i_1-1)
\Theta_J(k+l)
w_I
w_J
e_{I\!J}.
\]
By properties of the $\Theta$ function, we have
\begin{align*}
g(K)
&=\Theta_K(a+l)
w_K
e_K
=\Theta^-_{\overline{K}}(b-1)
w_K
e_K,
\quad\text{and}\\
X_G
&=
\sum_{K\vDash n}
(a-1)g(K)
-
\sum_{
2\le\abs{I}\le a-1,\
J\vDash n-\abs{I}}
(i_1-1)
\Theta^-_{\overline{J}}(b-1)
w_I
w_J
e_{I\!J}.
\end{align*}
For convenience, let
\[
\mathcal V_n=\{I\vDash n\colon i_1\ge 1,\
i_2,i_3,\dotsm\ge 2\}
\quad\text{and}\quad
\mathcal W_n=\{I\vDash n\colon i_1,i_2,\dotsm\ge 2\}.
\]
Since $(i_1-1)w_Iw_J=i_1w_{J\!I}$, we can rearrange $X_G$ as $X_G=X_1+X_2+X_3$, where
\begin{equation}\label{pf:CPC:def:Yu}
X_u
=
\sum_{
K\in\mathcal K_u}
(a-1)
g(K)
-
\sum_{
(I,J)\in\mathcal H_u}
\Theta^-_{\overline{J}}(b-1)
i_1
w_{J\!I}
e_{I\!J},
\end{equation}
for $u=1,2,3$, with
\begin{align*}
\mathcal K_1
&=\{
K\in\mathcal V_n\colon
k_1\le l+1
\},\\
\mathcal K_2
&=\{
K\in\mathcal W_n\colon
l+2\le k_1\le a+l-1
\},\\
\mathcal K_3
&=\{
K\in\mathcal W_n\colon
k_1\ge a+l
\},\\
\mathcal H_1
&=
\{
(I,J)\colon 
I\ne \emptyset,\
J\ne \emptyset,\
J\!I\in\mathcal V_n,\ \
\abs{I}\le a-1,\
j_1\le l+1\},\\
\mathcal H_2
&=
\{
(I,J)\colon 
I\ne \emptyset,\
J\ne \emptyset,\
J\!I\in\mathcal W_n,\
\abs{I}\le a-1,\
l+2\le j_1\le \abs{J}-b+1
\},\\
\mathcal H_3
&=
\{
(I,J)\colon 
I\ne \emptyset,\
J\ne \emptyset,\
J\!I\in\mathcal W_n,\
\abs{I}\le a-1,\
j_1\ge \abs{J}-b+2
\}.
\end{align*}
We observe that $\mathcal H_3=\mathcal H$, since
\[
j_1\ge\abs{J}-b+2
\iff
\abs{I}\ge n-b+2-j_1=a+l+1-j_1.
\]

For any $u\in \{1,2\}$ and any $(I,J)\in\mathcal H_u$, we consider the composition
\[
f(I,J)
=j_1\cdotp
I\cdotp
(J\backslash j_1).
\]
Then $f(I,J)\in \mathcal K_u$. It is immediate that $w_{J\!I}=w_{f(I,J)}$. Since $\abs{J\backslash j_1}\ge b-1$, we have
\[
\Theta^-_{\overline{f(I,J)}}(b-1)
=\Theta^-_{\overline{J}}(b-1).
\]
It follows that
\begin{equation}\label{pf:CPC:summand.-}
\Theta^-_{\overline{J}}(b-1)
i_1
w_{J\!I}
e_{I\!J}
=i_1
g(f(I,J)).
\end{equation}
Now we process $X_1$, $X_2$ and $X_3$ respectively.

First, we claim that
\begin{equation}\label{pf:CPC:Y1}
X_1
=\sum_{K\in \mathcal K_1}
s_1(K)g(K),
\quad\text{where }
s_1(K)
=\Theta_K^-(k_1+a-1).
\end{equation}
Let $K\in\mathcal K_1$, and write $\sigma_K^-(k_1+a-1)=k_1+\dots+k_{s+1}$ with $s\ge 0$. If $s\ge 1$, then $\mathcal H_1\cap f^{-1}(K)$ consists of the $s$ pairs obtained by taking $I=k_2\dotsm k_{r+1}$ for $1\le r\le s$; in every case  $i_1=k_2$. If $s=0$, this fibre is empty. Thus, by \cref{pf:CPC:summand.-},
\begin{equation}\label{pf:CPC:Y1'}
X_1
=\sum_{K\in\mathcal K_{1}}
(a-1-sk_2)
g(K).
\end{equation}
Let
\[
\mathcal K_1^0
=\{K\in\mathcal K_1\colon \sigma_{K\backslash k_1}^-(a-1)=0\}.
\]
For $K\in\mathcal K_1^0$, one has $s=0$ and $s_1(K)=a-1$, so these terms already have the form required in \cref{pf:CPC:Y1}.

In order to show \cref{pf:CPC:Y1} for the remaining terms, we consider the conjugacy equivalence relation on compositions: $P$ and~$Q$ are \emph{conjugate} if and only if $P=uv$ and $Q=vu$ for some compositions $u$ and $v$. Let $\mathcal W_n'$ be the set of conjugacy classes of $\mathcal W_n$. Then
\[
\mathcal K_{1}\backslash\mathcal K_1^0
=\bigsqcup_{k_1=1}^{l+1}
\bigsqcup_{d=2}^{a-1}
\bigsqcup_{C\in \mathcal W_d'}
\bigsqcup_{T\vDash n-k_1-d}
\mathcal K_{1}(k_1,d,T,C),
\]
where
\[
\mathcal K_{1}(k_1,d,T,C)
=\{K=k_1PT\in\mathcal V_n\colon
\sigma_{K\backslash k_1}^-(a-1)=d,\
P\in C\}.
\]
Let $P$ have length $s$ and let its cyclic conjugacy class $C$ have size $l'$. The cyclic orbit size $l'$ divides $s$, and $P=B^{s/l'}$ for a block $B$ of length $l'$. Thus, for any $K\in \mathcal K_1(k_1,d,T,C)$, the class $C$ can be written as
\[
C=\{(k_2k_3\dotsm k_{l'+1})^{s/l'},\
(k_3k_4\dotsm k_{l'+1}k_2)^{s/l'},\
\dots,\
(k_{l'+1}k_2k_3\dotsm k_{l'})^{s/l'}
\},
\]
and every composition in this class has the same total size
\begin{equation}\label{pf:CPC:Y1.k2sum}
(k_2+\dots+k_{l'+1})\cdotp s/l'
=\sigma_{K\backslash k_1}^-(a-1)
=a-1-s_1(K).
\end{equation}
On the other hand, we observe that $g(K)$ is invariant for $K\in\mathcal K_1(k_1,d,T,C)$. In fact,
\begin{enumerate}
\item
the factor $\Theta_{\overline K}^-(b-1)$ is invariant since
\[
\Theta_{\overline K}^-(b-1)
=\Theta_{\overline{k_1PT}}^-(b-1)
=\Theta_{\overline{T}\,\overline{P}\,k_1}^-(b-1)
=\Theta_{\overline T}^-(b-1)
\]
by the premise $\abs{T}=n-k_1-d\ge n-(l+1)-(a-1)=b-1$, and
\item
the factor $w_K e_K$ is invariant since the first part of $K$ is fixed to be $k_1$ and the other parts permute when $K$ runs over $\mathcal K_1(k_1,d,T,C)$.
\end{enumerate}
Now, for any $H\in\mathcal K_1(k_1,d,T,C)$, by \cref{pf:CPC:Y1.k2sum},
\begin{align*}
\frac{1}{\abs{C}}
\sum_{K\in\mathcal K_{1}(k_1,d,T,C)}
(a-1-sk_2)
g(K)
&=
\brk3{
(a-1)
-\frac{s}{\abs{C}}\brk1{
h_2+\dots+h_{\abs{C}+1}}}
g(H)
=
s_1(H)
g(H).
\end{align*}
Since $g$ is invariant in $\mathcal K_1(k_1,d,T,C)$, one may regard this identity as that every composition $K$ in $\mathcal K_{1}(k_1,d,T,C)$ contributes the term $s_1(K)g(K)$ on average. Together with the terms indexed by $\mathcal K_1^0$, \cref{pf:CPC:Y1'} yields \cref{pf:CPC:Y1}.

Second, we claim that
\begin{equation}\label{pf:CPC:Y2}
X_2
=\sum_{K\in \mathcal K_2}
s_2(K)
g(K),
\quad
\text{where }
s_2(K)
=
\Theta_K^-(a+l)+k_1-l-1.
\end{equation}
If $k_1=a+l-1$, then $s_2(K)=a-1$ and the summand in \cref{pf:CPC:Y2} coincides with the first summand in \cref{pf:CPC:def:Yu}. Thus we only need to focus on the set
\[
\mathcal K_2'
=\{K\in\mathcal W_n\colon
l+2\le k_1\le a+l-2\}
\]
for the positive sum of \cref{pf:CPC:def:Yu}.

We will show \cref{pf:CPC:Y2} in a way that is slightly different from that used for \cref{pf:CPC:Y1}. Let $K\in\mathcal K_2'$, and write
\[
\sigma_{K}^-(a+l)
=k_1+\dots+k_{t+1}
\]
with $t\ge 0$. If $t\ge 1$, then $\mathcal H_2\cap f^{-1}(K)$ consists of the $t$ pairs obtained by taking $I=k_2\dotsm k_{r+1}$ for $1\le r\le t$; if $t=0$, this fibre is empty. It follows that
\begin{equation}\label{pf:CPC:Y2.2}
X_2
=\sum_{K\in \mathcal K_2,\
k_1=a+l-1}
s_2(K)
g(K)
+
\sum_{
K\in\mathcal K_2'}
(a-1-tk_2)
g(K).
\end{equation}
Let
\[
(\mathcal K_2')^0
=\{K\in\mathcal K_2'\colon
\sigma_{K\backslash k_1}^-(a+l-k_1)=0\}.
\]
For $K\in(\mathcal K_2')^0$, one has $t=0$ and $s_2(K)=a-1$, so the corresponding term in \cref{pf:CPC:Y2.2} already has the required form. For the remaining terms, we have
\[
\mathcal K_2'\backslash(\mathcal K_2')^0
=\bigsqcup_{k_1=l+2}^{a+l-2}
\bigsqcup_{d=2}^{a+l-k_1}
\bigsqcup_{C\in \mathcal W_d'}
\bigsqcup_{
T\vDash n-k_1-d}
\mathcal K_2(k_1,d,T,C),
\]
where
\[
\mathcal K_2(k_1,d,T,C)
=\{K=k_1PT\in\mathcal W_n\colon
\sigma_{K\backslash k_1}^-(a+l-k_1)=d,\
P\in C\}.
\]
For any composition $K=k_1PT\in \mathcal K_2(k_1,d,T,C)$, the class $C$ can be written as
\[
C=\{(k_2k_3\dotsm k_{l'+1})^{t/l'},\
(k_3k_4\dotsm k_{l'+1}k_2)^{t/l'},\
\dots,\
(k_{l'+1}k_2k_3\dotsm k_{l'})^{t/l'}
\},
\]
where $l'=\abs{C}$, and every composition in this class has the same total size
\[
(k_2+\dots+k_{l'+1})\cdotp t/l'
=\sigma_K^-(a+l)-k_1
=a+l-\Theta_K^-(a+l)-k_1
=a-1-s_2(K).
\]
The function $g$ is invariant on $\mathcal K_2(k_1,d,T,C)$. Indeed, $d\le a+l-k_1$ gives $\abs{T}=n-k_1-d\ge b-1$, so the factor $\Theta^-_{\overline K}(b-1)$ depends only on $T$; the factor $w_Ke_K$ is unchanged under cyclic rotations of the block $P$. Therefore, for any $H\in\mathcal K_2(k_1,d,T,C)$,
\[
\frac{1}{\abs{C}}
\sum_{K\in\mathcal K_2(k_1,d,T,C)}
(a-1-tk_2)
g(K)
=
\brk3{
a-1
-\frac{t}{\abs{C}}
\brk1{h_2+\dots+h_{\abs{C}+1}}
}
g(H)
=
s_2(H)
g(H).
\]
Together with the terms indexed by $(\mathcal K_2')^0$ and by $k_1=a+l-1$, \cref{pf:CPC:Y2.2} yields \cref{pf:CPC:Y2}.

Now we deal with $X_3$. Let us keep in mind that $\mathcal H_3=\mathcal H$. Recall from \cref{pf:CPC:def:Yu} that $X_3$ has a positive sum
\[
X_3^+
=
\sum_{
K\in\mathcal K_3
}
(a-1)
g(K)
\]
and a negative sum
\[
-X_3^-
=
-
\sum_{
(I,J)\in\mathcal H}
\Theta^-_{\overline{J}}(b-1)
i_1
w_{J\!I}
e_{I\!J}.
\]
We shall produce a positive $e_I$-expansion for $X_G-X_3^+=X_1+X_2-X_3^-$.

First, we rewrite $-X_3^-$ in terms of $g(I\!J)$. We claim that
\begin{equation}\label{pf:def:-X3'}
-X_3^-
=-
\sum_{
(I,J)\in\mathcal H
}
i_1 s_3(I,J)
g(I\!J),
\quad\text{where }
s_3(I,J)
=\frac{j_1}{i_1}
\cdotp
\frac{i_1-1}{j_1-1}.
\end{equation}
In fact, for any $(I,J)\in\mathcal H$,
\[
g(I\!J)
=\Theta_{I\!J}(a+l)w_{I\!J}e_{I\!J}
=\Theta_{\overline{I\!J}}^-(b-1)w_{J\!I}e_{I\!J}/s_3(I,J),
\]
in which
\begin{equation}\label{pf:CPC:Y1.Theta2}
\Theta_{\overline{J}}^-(b-1)
=b-1-\abs{J\backslash j_1}
=\Theta_{\overline{I\!J}}^-(b-1),
\end{equation}
since $\abs{J\backslash j_1}<b-1<b+l=n-(a-1)\le n-\abs{I}=\abs{J}$. This proves the claim.

Now, the function $-X_3^-$ is expressed in terms of $g(I\!J)$. In order to merge it with $X_1+X_2$, in view of \cref{pf:CPC:Y1,pf:CPC:Y2}, we always consider
\[
K=I\!J
\quad\text{for all $(I,J)\in\mathcal H$}.
\]
For statement convenience, we define the concatenation map
\begin{equation}\label{def:concatenation}
h\colon \mathcal H\to
\mathcal W_n,
\quad\text{namely }
(I,J)\mapsto I\!J.
\end{equation}
We claim that $h$ is injective. In fact, suppose that $h(I,J)=h(P,Q)$ for some $(P,Q)\in\mathcal H$. If $I=P$, then $J=Q$ and we are done. Assume that $I\ne P$. Without loss of generality, we can suppose that $P=Ij_1\alpha$, for some composition $\alpha$ which might be empty. For any $(I,J)\in\mathcal H_3$, since
\begin{equation}\creflabel[ineq]{pf:CPC:Y3.I+j1}
\abs{I}
+j_1
=n-\abs{J\backslash j_1}
\ge n-(b-2)
=a+l+1,
\end{equation}
we find
\begin{equation}\label[ineq]{pf:CPC:Y3.j1>=l+2}
j_1
\ge a+l+1-\abs{I}
\ge l+2.
\end{equation}
By \cref{pf:CPC:Y3.I+j1}, we infer that
\[
a+l+1
\le\abs{I}+j_1
\le\abs{P}
\le a-1,
\]
which is absurd. This proves the injectivity of $h$.

We proceed according to the decomposition $\mathcal H=\mathcal H'\sqcup\mathcal H^-\sqcup\mathcal H^+$, where
\begin{align*}
\mathcal H'
&=\{
(I,J)\in\mathcal H\colon
2\le i_1\le l+1<j_1\}
=\{
(I,J)\in\mathcal H\colon
i_1\le l+1\},\\
\mathcal H^-
&=\{
(I,J)\in\mathcal H\colon
l+2\le j_1<i_1\}
=\{(I,J)\in\mathcal H\colon i_1>j_1\},
\\
\mathcal H^+
&=\{
(I,J)\in\mathcal H\colon
l+2\le i_1\le j_1\}.
\end{align*}
It is routine to check that
\[
h(\mathcal H')\subseteq\mathcal K_1
\quad\text{and}\quad
h(\mathcal H^-\cup\mathcal H^+)\subseteq\mathcal K_2.
\]
Since $h$ is injective, by \cref{pf:CPC:Y1,pf:CPC:Y2,pf:def:-X3'}, we obtain
\begin{equation}\label{pf:X1+X2-X3prime}
X_G
=
\sum_{K\in\mathcal K_1\backslash h(\mathcal H')}
s_1(K)g(K)
+
\sum_{K\in\mathcal K_2\backslash h(\mathcal H^-\cup\mathcal H^+)}
s_2(K)g(K)
+
X_3^+
+
Y_1
+
Y_2,
\end{equation}
where
\begin{align*}
Y_1
&=
\sum_{
(I,J)\in\mathcal H'
}
\brk1{
s_1(I\!J)
-
i_1 s_3(I,J)
}
g(I\!J),
\quad\text{and}\\
Y_2
&=
\sum_{
(I,J)\in\mathcal H^-\sqcup\mathcal H^+
}
\brk1{
s_2(I\!J)
-
i_1 s_3(I,J)
}
g(I\!J).
\end{align*}
We shall produce a positive $e_I$-expansion for $Y_1$ and $Y_2$, respectively.
\begin{enumerate}
\item
Let $(I,J)\in\mathcal H'$. By \cref{pf:CPC:Y3.I+j1}, we have $\abs{I}\le a-1<a-1+i_1\le a+l<\abs{I}+j_1$. Thus
\begin{equation}\label{pf:CPC:Y3.s1}
s_1(I\!J)
=\Theta_{I\!J}^-(i_1+a-1)
=i_1+a-1-\abs{I}.
\end{equation}
It follows that $Y_1=\sum_{(I,J)\in\mathcal H'}t_1(I,J)g(I\!J)$, where
\begin{align}
\label{pf:t1:kayak}
t_1(I,J)
=
s_1(I\!J)
-
i_1
s_3(I,J)
=
a-1-\abs{I}
+
\frac{j_1-i_1}
{j_1-1}.
\end{align}
Since $j_1>i_1$, we conclude that $t_1(I,J)>0$.
\item
Let $(I,J)\in\mathcal H^-\cup\mathcal H^+$. By \cref{pf:CPC:Y3.I+j1}, we have $\abs{I}<a+l<\abs{I}+j_1$ and $\Theta_{I\!J}^-(a+l)=a+l-\abs{I}$. It follows that
\[
s_2(I\!J)
=\Theta_{I\!J}^-(a+l)+i_1-l-1
=a-\abs{I}+i_1-1,
\]
which shares the same expression on the right side of \cref{pf:CPC:Y3.s1} for $s_1(I\!J)$. Therefore, by \cref{pf:t1:kayak},
\[
s_2(I\!J)
-
i_1
s_3(I,J)
=
t_1(I,J).
\]
It is nonnegative for $(I,J)\in\mathcal H^+$. We need to deal with $(I,J)\in\mathcal H^-$.

Consider the involution $\phi$ defined for pairs $(I,J)$ of nonempty compositions by
\[
\phi(I,J)=(P,Q),
\quad
\text{where
$P=j_1(I\backslash i_1)$
and
$Q=i_1(J\backslash j_1)$}.
\]
We claim that $\phi(\mathcal H^-)\subseteq \mathcal H^+$. In fact, for any $(I,J)\in\mathcal H^-$, we have $i_1>j_1\ge l+2$; moreover,
\begin{itemize}
\item
every part of $PQ$ is at least $2$, as are the parts of $I\!J$,
\item
$q_1=i_1>j_1=p_1$ and $p_1=j_1\ge l+2$,
\item
$\abs{Q}-q_1=\abs{J}-j_1\le b-2$, and
\item
$2\le j_1\le\abs{P}=\abs{I}-i_1+j_1<\abs{I}\le a-1$.
\end{itemize}
This proves the claim. On the other hand, recall from \cref{pf:CPC:Y1.Theta2} that
\[
\Theta_{\overline{I\!J}}^-(b-1)
=b-1-\abs{J\backslash j_1}
=b-1-\abs{Q\backslash q_1}
=\Theta_{\overline{PQ}}^-(b-1).
\]
Since $w_{PQ}=s_3(I,J)w_{I\!J}$, we can compute that
\[
t_1(I,J)g(I\!J)
+t_1(P,Q)g(PQ)
=t_2(I,J)g(I\!J),
\]
where
\begin{align*}
t_2(I,J)
&=
\brk3{
a-1-\abs{I}
+
\frac{j_1-i_1}
{j_1-1}
}
+
\brk3{
a-1-\brk1{
\abs{I}-i_1+j_1}
+
\frac{i_1-j_1}
{i_1-1}
}
s_3(I,J)\\
&=
(i_1-j_1)
+
\brk1{
a-1-\abs{I}}
\brk1{
1+s_3(I,J)
}
>0.
\end{align*}
Hence we obtain a positive $e_I$-expansion
\[
Y_2
=
\sum_{
(I,J)\in
\mathcal H^+
\backslash \phi(\mathcal H^-)
}
t_1(I,J)
g(I\!J)
+
\sum_{
(I,J)\in\mathcal H^-
}
t_2(I,J)
g(I\!J).
\]
\end{enumerate}
Before combining the right side of \cref{pf:X1+X2-X3prime}, we simplify the difference sets
\[
\mathcal K_1\backslash h(\mathcal H'),
\quad
\mathcal K_2\backslash h(\mathcal H^-\cup\mathcal H^+),
\quad\text{and}\quad
\mathcal H^+
\backslash \phi(\mathcal H^-).
\]
\begin{enumerate}
\item
Recall that $\mathcal K_1=\{K\in\mathcal V_n\colon k_1\le l+1\}$ and
\[
\mathcal H'
=\{(I,J)\colon 
J\!I\in\mathcal W_n,\
I\ne\emptyset,\
J\ne\emptyset,\
a+l-j_1+1
\le \abs{I}
\le a-1,\
i_1\le l+1
\}.
\]
Let $K\in\mathcal K_1\backslash h(\mathcal H')$. Suppose that $k_1\ge 2$. Then $K\in\mathcal W_n$ and $k_1\le l+1$. We infer that if $K=I\!J$ for some $I,J\ne\emptyset$ with $\abs{I}\le a-1$, then $\abs{I}\le a+l-j_1$, i.e., $\abs{J\backslash j_1}\ge b-1$. This is equivalent to saying that if $\abs{I}=\sigma_K^-(a-1)$ and $I=k_1\dotsm k_\alpha$ for some $\alpha\ge 1$, then $j_1=k_{\alpha+1}$ and
\[
b-1\le \abs{J\backslash j_1}=n-\abs{I}-j_1=n-\abs{k_1\dotsm k_{\alpha+1}},
\]
i.e., $a+l\ge \abs{k_1\dotsm k_{\alpha+1}}=\sigma_K(a)$, or equivalently, $\Theta_K(a)\le l$. Consequently,
\begin{equation}\label{pf:K1-hH3}
\mathcal K_1\backslash h(\mathcal H')
=
\{
K\in\mathcal K_1
\colon 
k_1=1
\text{ or }
\Theta_K(a)\le l
\}.
\end{equation}
\item
Along the same lines, we deduce that
\begin{equation}\label{pf:K2-hH3}
\mathcal K_2\backslash h(\mathcal H^-\cup\mathcal H^+)
=
\{
K\in\mathcal K_2
\colon
\Theta_K(a)\le l
\}.
\end{equation}
\item
Note that
\begin{align*}
\mathcal H^-
&=
\{(I,J)
\colon
J\!I\in\mathcal W_n,\
2\le\abs{I}\le a-1,\
j_1\ge \abs{J}-b+2,\
l+2\le j_1<i_1\}
\\
&=
\{(I,J)
\colon
I\!J\in\mathcal W_n,\
\max(l+2,\,
\abs{J}-b+2)
\le j_1
<i_1
\le
\abs{I}
\le a-1
\}.
\end{align*}
If $(I,J)\in\mathcal H^-$, then neither $I$ nor $J$ is empty. If $(P,Q)=\phi(I,J)$, then neither $P$ nor $Q$ is empty, and $p_1=j_1$, $q_1=i_1$, $\abs{P}=\abs{I}+p_1-q_1$, and $\abs{Q}=\abs{J}-p_1+q_1$. It follows that
\begin{align*}
\phi(\mathcal H^-)
&=
\left\{(P,Q)\colon
\begin{gathered}
P Q\in\mathcal W_n,\quad P\ne\emptyset,\quad Q\ne\emptyset,\\
\max(l+2,\,\abs{Q}+p_1-q_1-b+2)\le p_1<q_1,\\
q_1\le \abs{P}-p_1+q_1\le a-1
\end{gathered}
\right\}
\\
&=
\left\{(P,Q)\colon
\begin{gathered}
P Q\in\mathcal W_n,\quad P\ne\emptyset,\quad Q\ne\emptyset,\\
\abs{Q}-b+2\le q_1,\quad l+2\le p_1<q_1,\\
\abs{P}-p_1+q_1\le a-1
\end{gathered}
\right\}
\\
&=
\left\{(I,J)\colon
\begin{gathered}
I\!J\in\mathcal W_n,\quad I\ne\emptyset,\quad J\ne\emptyset,\\
\abs{J}-b+2\le j_1,\quad l+2\le i_1<j_1,\\
\abs{I}\le a-1+i_1-j_1
\end{gathered}
\right\}.
\end{align*}
Since
\begin{align*}
\mathcal H^+
&=
\{(I,J)
\colon
J\!I\in\mathcal W_n,\
2\le\abs{I}\le a-1,\
j_1\ge \abs{J}-b+2,\
j_1\ge i_1\ge l+2\}
\\
&=
\{(I,J)
\colon
J\!I\in\mathcal W_n,\
I\ne \emptyset,\
J\ne \emptyset,\
\abs{J}-b+2\le j_1,\
l+2\le i_1\le j_1,\
\abs{I}\le a-1
\},
\end{align*}
we deduce that
\begin{align*}
\mathcal H^+
\backslash \phi(\mathcal H^-)
&=
\left\{(I,J)\colon
\begin{gathered}
J\!I\in\mathcal W_n,\quad \abs{I}\le a-1,\\
j_1=i_1\ge \max(l+2,\,\abs{J}-b+2)
\end{gathered}
\right\}
\\
&\quad\cup
\left\{(I,J)\colon
\begin{gathered}
J\!I\in\mathcal W_n,\quad I\ne\emptyset,\quad J\ne\emptyset,\\
\abs{J}-b+2\le j_1,\quad l+2\le i_1<j_1,\\
a+i_1-j_1\le \abs{I}\le a-1
\end{gathered}
\right\}
\\
&=
\{(I,J)\in\mathcal H^+\colon
i_1=j_1
\text{ or }
\abs{I}-i_1+j_1\ge a
\}.
\end{align*}
Since $\mathcal H'\cap \mathcal H^+=\emptyset$ and $\phi(\mathcal H^-)\subseteq \mathcal H^+$, we obtain
\begin{align*}
\mathcal H'
\cup
\brk1{
\mathcal H^+
\backslash \phi(\mathcal H^-)
}
&=\mathcal H'
\cup
\mathcal H^+
\backslash \phi(\mathcal H^-)
\\
&=
\{(I,J)\in\mathcal H\colon
i_1\le l+1
\text{ or }
i_1=j_1
\text{ or }
(i_1\ge l+2 \text{ and }\abs{I}\ge a+i_1-j_1)
\}.
\end{align*}
\end{enumerate}
Now, by \cref{pf:X1+X2-X3prime,pf:K1-hH3,pf:K2-hH3},
\begin{multline*}
X_G
=
\sum_{
\substack{
K\vDash n,\
\text{either } k_1=1\\
\text{or ($2\le k_1\le l+1$ and $\Theta_K(a)\le l$)}
}}
s_1(K)g(K)
+
\sum_{
K\in\mathcal K_2,\
\Theta_K(a)\le l
}
s_2(K)g(K)
+
\sum_{
K\in\mathcal K_3
}
(a-1)
g(K)
\\
+
\sum_{
\substack{
(I,J)
\in \mathcal H\\
i_1\le l+1
\text{ or }
i_1=j_1
\text{ or }
(i_1\ge l+2 \text{ and }\abs{I}\ge a+i_1-j_1)
}}
t_1(I,J)
g(I\!J)
+
\sum_{
(I,J)\in\mathcal H,\
i_1>j_1
}
t_2(I,J)
g(I\!J).
\end{multline*}
Note that the terms containing $e_1$ all come from the first sum, and they form the sum
\[
\sum_{K\vDash n,\ k_1=1}
s_1(K)g(K)
=
\sum_{K\vDash n,\ k_1=1}
\Theta_K^-(a)g(K).
\]
On the other hand, when $\Theta_K(a)=l$, we have $g(K)=\Theta_K(a+l)w_Ke_K=0$ and the condition $\Theta_K(a)\le l$ can be replaced equivalently with $\Theta_K(a)\le l-1$. Writing out the first sum for $k_1=1$ and $2\le k_1\le l+1$, and the functions $s_1$, $s_2$, $t_1$ and $t_2$, we obtain the desired formula.
\end{proof}

The kayak paddle $P^0(C_a,C_b)$ resembles the symbol $\infty$; we therefore call it the \emph{$(a,b)$-infinity graph} and denote it by $\infty_{ab}$. It is obtained by identifying a vertex of $C_a$ with a vertex of $C_b$. For any composition $I\vDash n$ and any number $1\le a\le n-1$, we denote
\[
I(a)=\Theta_I(a)+\Theta_I^-(a),
\] 
or equivalently, $I(a)=\sigma_I(a)-\sigma_I^-(a)$. In other words, if $a$ is a partial sum of $I$, then $I(a)=0$; otherwise, the value of $I(a)$ equals the part $i_k$ such that
\[
\abs{i_1\dotsm i_{k-1}}
<a
<\abs{i_1\dotsm i_k}.
\]
For example, for the composition $I=83617$ we have $I(15)=6$ and $I(17)=0$.

\begin{corollary}[Infinity graphs]\label{cor:infinity}
Let $a,b\ge 3$ and $n=a+b-1$. Then the chromatic symmetric function of the infinity graph $\infty_{ab}$ has the following positive $e_I$-expansion:
\begin{align*}
X_{\infty_{ab}}
&=
\sum_{
\substack{
I\vDash n\\
i_1=1
}}
\Theta_I^-(a)g(I)
+
\sum_{
\substack{
I\vDash n,\
2\le i_1\le a-1
\\
i_1\le\Theta_I(a)
\textnormal{ or }
i_1=I(a)
}}
\brk3{
\Theta_I^-(a)-1
+
\frac{I(a)-i_1}
{I(a)-1}}
g(I)
\\
&\quad
+
\sum_{
\substack{
I\vDash n\\
3\le I(a)+1\le i_1\le a-1
}}
\brk4{
i_1-I(a)
+
\brk1{
\Theta_I^-(a)
-1}
\brk3{
1+
\frac{I(a)}{i_1}
\cdotp
\frac{i_1-1}{I(a)-1}
}}
g(I)
\\
&\quad+
\sum_{
I\vDash n, \
i_1\ge a
}
(a-1)
g(I),
\end{align*}
where $g(I)=\Theta_I(a)w_Ie_I$.
\end{corollary}
\begin{proof}
Let $G=\infty_{ab}$. Taking $l=0$ in \cref{thm:kayak}, we obtain
\[
X_G
=
\sum_{
K\vDash n,\
k_1=1
}
\Theta_K^-(a)g(K)
+
\sum_{
K\vDash n,\
k_1\ge a
}
(a-1)
g(K)
+A+B,
\]
where $g(K)=\Theta_K(a)w_Ke_K$,
\begin{align*}
A&=
\sum_{
(I,J)
\in \mathcal H,\
i_1=j_1 \text{ or } 
\abs{I}\ge a+i_1-j_1
}
\brk3{
a-1-\abs{I}
+
\frac{j_1-i_1}
{j_1-1}}
g(I\!J),
\\
B
&=
\sum_{
(I,J)\in\mathcal H,\
i_1>j_1
}
\brk4{
(i_1-j_1)
+
\brk1{
a-1-\abs{I}}
\brk3{
1+
\frac{j_1}{i_1}
\cdotp
\frac{i_1-1}{j_1-1}
}}
g(I\!J),
\quad\text{and}\\
\mathcal H
&=\{(I,J)\colon
I\!J\in\mathcal W_n,\
I\ne\emptyset,\
J\ne\emptyset,\
a-j_1+1
\le \abs{I}
\le a-1
\}.
\end{align*} 
We shall simplify $A$ and $B$ by reading the pairs $(I,J)$ from a single composition $I\!J$, which ranges over a set that will be derived below. Let $\mathcal K=\{K\in \mathcal W_n\colon k_1\le a-1,\ \Theta_K(a)\ge 1\}$.

First, we claim that $\mathcal K=h(\mathcal H)$, where $h$ is the concatenation map defined by \cref{def:concatenation}. In fact, it is clear that $h(\mathcal H)\subseteq\mathcal K$. Conversely, any composition $K\in\mathcal K$ can be decomposed uniquely as $I\!J$, where $\abs{I}=\sigma_K^-(a)$. Under this decomposition, we check $(I,J)\in\mathcal H$ as follows.
\begin{itemize}
\item
Since $k_1\le a-1$, we have $I\ne\emptyset$.
\item
The premise $\Theta_K(a)\ge1$ is equivalent to $\Theta_K^-(a)\ge1$, which implies $\abs{I}\le a-1$.
\item
Since $a+b-1=\abs{I}+\abs{J}$, $\abs{I}\le a-1$ and $b\ge 2$, we find $J\ne \emptyset$.
\item
Since $\abs{I}=\sigma_K^-(a)$ and $\Theta_K(a)\ge1$, we find $\abs{I}+j_1=\sigma_K(a)\ge a+1$.
\end{itemize}
This proves $\mathcal K\subseteq h(\mathcal H)$ and the claim follows.

Now, for $K\in\mathcal K$, one may read the prefix $I$ of $K$ such that $\abs{I}=\sigma_K^-(a)$. Consequently, $i_1=k_1$, $j_1=K(a)$, $a-\abs{I}=\Theta_K^-(a)$, and the condition $\abs{I}\ge a+i_1-j_1$ becomes
\[
k_1
=i_1
\le\abs{I}-a+j_1
=
-\Theta_K^-(a)
+K(a)
=\Theta_K(a).
\]
Therefore, for the pairs $(I,J)$ in $A$ and $B$, the concatenations~$I\!J$ run over the sets
\begin{equation}\label{pf:K'.AB}
\{
K\in\mathcal K
\colon
k_1\le\Theta_K(a)
\text{ or }
k_1=K(a)
\}
\quad\text{and}\quad
\{
K\in\mathcal K
\colon
k_1\ge K(a)+1
\},
\end{equation}
respectively. Since $\Theta_K(a)\le K(a)-1$, these sets are disjoint.

Finally, we claim that both occurrences of the set $\mathcal K$ in \eqref{pf:K'.AB} can be replaced with the set
\[
\{K\vDash n\colon 2\le k_1\le a-1,\ K(a)\ge 2\}.
\]
In fact, if $\Theta_K(a)=0$ or $k_j=1$ for some $j\ge 2$, then the factor~$g(K)$ vanishes and so do $A$ and~$B$. On the other hand, since the difference $j_1-1=K(a)-1$ appears in the denominators of $A$ and~$B$, we need to restrict $K(a)\ne 1$. This can be done by requiring $K(a)\ge 2$, since $K(a)=0$ only if $\Theta_K(a)=0$. Collecting this information, reordering the four sums in $X_G$ according to the value of~$k_1$, and substituting $K$ by the symbol $I$, we obtain the desired formula for $X_G$.
\end{proof}

We remark that the four sums in the formula of \cref{cor:infinity} are for distinct compositions $I$.

\section*{Acknowledgment}
This paper was completed while the second author was visiting Professor Jean-Yves Thibon at LIGM, Université Gustave Eiffel. The second author thanks Professor Thibon and LIGM for their hospitality.

\section*{Data availability statement}
Data sharing is not applicable because no datasets were generated or analyzed in this study.

\bibliography{../csf}

\end{document}